\RequirePackage[2020-02-02]{latexrelease}
\documentclass[11pt]{amsart}
\usepackage{amsfonts}
\usepackage{amsmath}
\usepackage{amsthm}
\usepackage{amssymb}
\usepackage{wasysym}
\usepackage[numbers]{natbib}
\usepackage[fit]{truncate}
\usepackage{hyperref}
\usepackage{theoremref}

\usepackage{mathtools}

\usepackage[dvipsnames]{xcolor}

\usepackage{tikz-cd}

\usepackage{thmtools}
\usepackage{thm-restate}

\usepackage{fmtcount}

\usepackage{enumitem}

\usepackage[UKenglish]{babel}
\usepackage[UKenglish]{isodate}

\usepackage{verbatim}
\usepackage[utf8]{inputenc}

\usepackage{dsfont}

\usepackage{etoolbox}
\makeatletter
\patchcmd{\@maketitle}
{\ifx\@empty\@dedicatory}
{\ifx\@empty\@date \else {\vskip3ex \centering\footnotesize\@date\par\vskip1ex}\fi
	\ifx\@empty\@dedicatory}
{}{}
\patchcmd{\@adminfootnotes}
{\ifx\@empty\@date\else \@footnotetext{\@setdate}\fi}
{}{}{}
\makeatother

\makeatletter
\def\namedlabel#1#2{\begingroup
	#2%
	\def\@currentlabel{#2}%
	\phantomsection\label{#1}\endgroup
}
\makeatother

\usepackage{catoptions}
\makeatletter

\def\Autoref#1{%
	\begingroup
	\edef\reserved@a{\cpttrimspaces{#1}}%
	\ifcsndefTF{r@#1}{%
		\xaftercsname{\expandafter\testreftype\@fourthoffive}
		{r@\reserved@a}.\\{#1}%
	}{%
		\ref{#1}%
	}%
	\endgroup
}
\def\testreftype#1.#2\\#3{%
	\ifcsndefTF{#1autorefname}{%
		\def\reserved@a##1##2\@nil{%
			\uppercase{\def\ref@name{##1}}%
			\csn@edef{#1autorefname}{\ref@name##2}%
			\autoref{#3}%
		}%
		\reserved@a#1\@nil
	}{%
		\autoref{#3}%
	}%
}
\makeatother

\numberwithin{equation}{section}

 \newtheoremstyle{dotless}{}{}{\itshape}{}{\bfseries}{}{ }{}

\theoremstyle{dotless}

\MakeRobust{\ref}

\makeatletter
\newcommand{\labeltext}[2]{%
	\@bsphack
	\csname phantomsection\endcsname 
	\def\@currentlabel{#1}{\label{#2}}%
	\@esphack
}
\makeatother

\makeatletter
\def\blfootnote{\gdef\@thefnmark{}\@footnotetext}
\makeatother

\newcommand{\N}[0]{\mathbb{N}}
\newcommand{\Z}[0]{\mathbb{Z}}
\newcommand{\Q}[0]{\mathbb{Q}}
\newcommand{\R}[0]{\mathbb{R}}

\newcommand{\supp}[0]{\mathrm{supp}}

\newcommand{\brackets}[1]{\left( #1 \right)}
\newcommand{\setbr}[1]{\left\{ #1 \right\}}
\newcommand{\pow}[1]{\!\left(\!\left( #1 \right)\!\right)}
\newcommand{\rpow}[1]{\!\left[\!\left[ #1 \right]\!\right]}

\newcommand{\ul}[1]{\underline{#1}}
\newcommand{\gen}[1]{\left\langle{#1}\right\rangle}
\newcommand\restr[2]{{
		\left.\kern-\nulldelimiterspace 
		#1
		\vphantom{\big|} 
		\right|_{#2}
}}

\newcommand{\x}{\underline{x}}

\DeclareMathOperator{\Char}{char}

\DeclareMathOperator{\ff}{ff}

\DeclareMathOperator{\Span}{span}

\newcommand{\bbF}{{\mathbb F}}

\newcommand{\cF}{{\mathcal F}}

\newcommand{\cP}{{\mathcal P}}

\newcommand{\C}{\mathbb C}

\newcommand{\one}{\mathds{1}}

\DeclareMathOperator{\id}{id}
\DeclareMathOperator{\Aut}{Aut}

\DeclareMathOperator{\Int}{Int}

\DeclareMathOperator{\lex}{lex}

\DeclareMathOperator{\Hom}{Hom}
\newcommand{\vAut}{v\text{-}\!\Aut}
\newcommand{\oAut}{\mathrm{o}\text{-}\!\Aut}

\DeclareMathOperator{\dom}{dom}
\DeclareMathOperator{\ot}{ot}

\newcommand{\Hk}[1]{\mathcal{H}^{#1}}

\DeclareMathOperator*\lc{LRS}
\DeclareMathOperator*\LRS{LRS}
\DeclareMathOperator*\fin{fin}

\renewcommand{\epsilon}{\varepsilon}
\renewcommand{\phi}{\varphi}
\renewcommand{\theta}{\vartheta}

\usepackage{theoremref}

\theoremstyle{definition}
\newtheorem{defn}{Definition}[section]
\newtheorem{definition}[defn]{Definition}

\newtheorem{notation}[defn]{Notation}
\theoremstyle{plain}
\newtheorem{lemma}[defn]{Lem\-ma}
\theoremstyle{plain}
\newtheorem*{lemma*}{Lem\-ma}
\theoremstyle{plain}

\newtheorem{proposition}[defn]{Proposition}
\newtheorem{fact}[defn]{Fact}

\newtheorem{question}[defn]{Question}
\newtheorem{corollary}[defn]{Corollary}
\theoremstyle{plain}
\theoremstyle{plain}
\newtheorem*{prop*}{Proposition}
\newtheorem{theorem}[defn]{Theorem}
\theoremstyle{plain}
\newtheorem*{teorema*}{Theorem}
\theoremstyle{plain}

\theoremstyle{plain}
\newtheorem*{cor*}{Corollary}
\theoremstyle{definition}

\newtheorem{example}[defn]{Example}
\theoremstyle{definition}

\theoremstyle{plain}

\theoremstyle{definition}

\theoremstyle{definition}

\newtheorem{construction}[defn]{Construction}
\newtheorem{remark}[defn]{Remark}
\theoremstyle{plain}
\theoremstyle{remark}

{\left\lbrace\begin{array}{@{}l@{}}}%
	{\end{array}\right.}

\newcommand{\vertiii}[1]{{\left\vert\kern-0.25ex\left\vert\kern-0.25ex\left\vert #1 
		\right\vert\kern-0.25ex\right\vert\kern-0.25ex\right\vert}}

\begin{document}

	\title[Generalised power series determined by linear recurrence]{Generalised power series determined by linear recurrence relations}
	
	\author[L.~S.~Krapp]{Lothar Sebastian Krapp}\thanks{Part of this work was carried out while all three authors were generously hosted by the Fields Institute, in Toronto, during the month of June 2022. We wish to thank Micka{\"e}l Matusinski and Victor Vinnikov, {as well as the anonymous referee}, for helpful discussions and comments.}
	\author[S.~Kuhlmann]{Salma Kuhlmann}
	\author[M.~Serra]{Michele Serra}
	
	\address{Institut für Interdisziplinäre Sprachevolutionswissenschaft, Universität Zürich, Affolternstrasse 56, 8050 Zürich, Switzerland \& Fachbereich Ma\-the\-ma\-tik und Statistik, Universität Konstanz, 78457 Konstanz, Germany}
	\email{sebastian.krapp2@uzh.ch}
	
	\address{Fachbereich Mathematik und Statistik, Universität Konstanz, 78457 Konstanz, Germany}
	\email{salma.kuhlmann@uni-konstanz.de}
	
	\address{Fakultät für Ma\-the\-ma\-tik, Technische Universität Dortmund, Vogelpoths\-weg 87, 44227 Dortmund, Germany}
	\email{{michele.serra@tu-dortmund.de}}
	
	\date{\today}
	
	\cleanlookdateon
	
	\maketitle
	
	\begin{abstract}
		
		In 1882, Kronecker established that a given univariate for\-mal Laurent series over a field  can be expressed as a fraction of two univariate polynomials if and only if the coefficients of the series satisfy a linear recurrence relation. We introduce the notion of \emph{generalised} linear recurrence relations for power series with exponents in an ar\-bi\-tra\-ry ordered abelian group, and generalise Kronecker's original result. In particular, we obtain criteria for determining whether a multivariate formal Laurent series lies in the fraction field of the corresponding polynomial ring. Moreover, we study distinguished algebraic substructures of a power series field, which are determined by generalised linear recurrence relations.
		In particular, we identify generalised linear recurrence relations that determine power series fields satisfying additional properties which are essential for the study of their automorphism groups.
		
	\end{abstract}
	
	
	\blfootnote{\textup{2020} \textit{Mathematics Subject Classification}: 13J05 (16W60 12J10 06F20)
		}

\section{Introduction}\label{sec:intro}

	For a field $k$ and an ordered abelian group $G$, a Hahn field is any subfield of the field of generalised power series $k\pow{G}$ (the maximal Hahn field) containing the minimal Hahn field $k(G)$, the fraction field of the group ring $k[G]$ (see Section \ref{sec:preliminaries} for details). Our original motivation for developing a notion of generalised linear recurrence relations stems from \cite{kuhlmann-serra-fields}, in which the second and third author introduce several so-called lifting properties (see Section \ref{sec:lifting}) in order to examine  the automorphism groups of Hahn fields. Indeed, we establish in this paper that Hahn fields whose elements satisfy specific generalised linear recurrence relations naturally exhibit these lifting properties (see \Autoref{thm:1lpcriterion} and \Autoref{thm:2lpcriterion}), which in turn leads to construction methods of such Hahn fields (see \Autoref{constr:1clpfields}).
	
	Besides considering lifting properties, we aim more globally to describe distinguished Hahn fields. In \cite{kks}, we focus on Hahn fields consisting of all series whose support lies in some family $\cF$ of well-ordered subsets of $G$, and we characterise substructures of Hahn fields via properties of $\cF$. Here, our description is focused on Hahn fields characterised by generalised linear recurrence relations involving both the support \emph{and} the coefficients of the series (see \Autoref{sec:hahnfields}).
	
	By building up the theory of generalised linear recurrence relations, we discovered that they represent a powerful tool in the study of Hahn fields and properties of their elements; due to Kronecker~\cite{kronecker}, an element of the field of formal Laurent series $k\pow{t}=k\pow{\Z}$ lies in $k(t)$ if and only if its coefficients satisfy a (finite) linear recurrence relation with constant coefficients (i.e.\ coefficients in $k$; see \Autoref{prop:laurentseries}).
	We generalise this result to any exponent group $G$ rather than $\Z$, i.e.\ we show that a generalised power series in $k\pow{G}$ already lies in $k(G)$ if and only if its coefficients satisfy a (finite) \emph{generalised} linear recurrence relation with constant coefficients (see \Autoref{lem:coefficientskG}). By a careful choice of the ordered exponent group $G$, \Autoref{lem:coefficientskG} can be applied to the ring of multivariate formal power series $k\rpow{\x}=k\rpow{x_1,\ldots,x_n}$ and its fraction field $k\pow{\x}$ (see \Autoref{theo-multivar} and \Autoref{prop:k[[x]]}).
	Moreover, \Autoref{constr:pq} presents a procedure to express a power series in $k\pow{\x}$ satisfying a linear recurrence relation as a fraction of two polynomials in $k[\x]$.
	Another characterisation of $k(t)$ within $k\pow{t}$ also presented by Kronecker~\cite{kronecker} utilised Hankel matrices: namely, a formal Laurent series in $k\pow{t}$ already lies in $k(t)$ if and only if its associated Hankel matrix has finite rank (see \Autoref{fact:kroneckerhankel}). We also obtain a generalisation of this characterisation in \Autoref{sec:gen-Hank-Matrices} (see \Autoref{fact:kroneckerhankelgen}).
	
	The structure of this paper is as follows. In \Autoref{sec:preliminaries} we set up the framework on ordered abelian groups and Hahn fields that we work in.
	For the convenience of the reader, we briefly review the original results of Kron\-ec\-ker~\cite{kronecker}. 
	In \Autoref{sec:gen-lin-rec}, we introduce and illustrate by examples our main concept of (generalised) linear recurrence sequences (\Autoref{def:lrs}) and sets determined by them (\Autoref{def:determined}). Moreover, we give a series of auxiliary lemmas leading to the proof of the main theorem of this section (\Autoref{prop:field}), which gives sufficient conditions on a set of linear recurrence sequences to determine a field.
	Thus, \Autoref{prop:field} singles out certain sets of linear recurrence sequences. We make this precise by introducing in \Autoref{sec:hahnfields} the key concept of $\cF$-sequences, which are sets of linear recurrence sequences whose domain lies in a family $\cF$ of well-ordered subsets of $G$ (\Autoref{def:fseq}). 
	In \Autoref{sec:gen-Hank-Matrices}, we introduce the concept of Hankel matrices of linear recurrence sequences. As an application, we recover a generalisation of the classical result that the exponential function is not rational (\Autoref{ex:k(G)-not-k-hull0}).
	In \Autoref{sec:rayner}, we relate Hahn fields determined by $\cF$-sequences to Rayner fields (which we introduced in \cite{kks}).
	The key \Autoref{rmk:k((x))props} of \Autoref{sec:multivariate} establishes that $k\pow{\x}$ is determined by linear recurrence relations and, for $n\geq 2$, is neither the minimal nor the maximal Hahn field.
	In \Autoref{sec:lifting1} and \Autoref{sec:lifting2} of Section \ref{sec:lifting}, we examine the canonical first and, respectively, the canonical second lifting property. While every Hahn field determined by $\cF$-sequences has the canonical \emph{second} lifting property (\Autoref{thm:canonicallifting2}), we show that some Rayner fields (\Autoref{ex:raynernolift}) and  $k\pow{{\x}}$ (\Autoref{ex:810}) are Hahn fields determined by $\cF$-sequences \emph{without} the canonical \emph{first} lifting property. We do not know whether $k\pow{\x}$ is also a Rayner field (\Autoref{qu:rayner}).
	
	Finally, we point out that our study of algebraic properties of formal power series in terms of their coefficients relates to several other examinations of formal power series and related concepts in the literature;
	determining algebraic properties of Puiseux series in terms of their coefficients was the subject of study of \cite{hickel} in the univariate case and of \cite{hickel2} in the multivariate case. Several other works give criteria on the rationality of a power series in terms of Hankel matrices: For instance, \cite{power} deals with the multivariate case over the coefficient field $\C$, \cite{duchamp} studies the case of non-commutative
	formal power series, which also connect to theoretical computer science via automata (see also \cite{reutenauer}), and Hankel operators of finite rank play an important role in the study of moment problems (see \cite{mourrain} for an extensive list of references).
	
\section{Preliminaries}\label{sec:preliminaries}

	Throughout this work, let $k$ be a field, let $G$ be a (totally) ordered abelian group, and let $\cF$ be a family of well-ordered subsets of $G$.
	
\subsection{Notation and terminology}

	We denote by $\N$ the set of natural numbers with $0$, by $G^{>0}$ the set of positive elements of $G$ and by $G^{\geq 0}$ the set of non-negative elements of $G$. 
	The power set of a set $B$ is denoted by $\mathcal{P}(B)$.
	The fraction field of an integral domain $R$ is denoted by $\ff(R)$.
	For any sets $A\subseteq G$ and $C\subseteq \Z$, we denote by $\Span_C(A)$ the subset of $G$ consisting of all finite linear combinations of the form $\sum_{i=1}^n c_ia_i$, where $n\in \N$, $c_i\in C$ and $a_i\in A$. Thus, the subgroup of $G$ generated by $A$ is given by $\Span_\Z(A)$. We use standard abbreviations like $\Span_C(A,h)$ for $\Span_C(A\cup\{h\})$ and $\Span_C(h,h')$ for $\Span_C(\{h,h'\})$, where $A\subseteq G$ and  $h,h'\in G$.
	
	For a partial function $f\colon G\rightharpoonup k$, the \textbf{support} of $f$ is given by $\supp(f)={\dom(f)\setminus f^{-1}(0)}=\{g\in \dom(f)\mid f(g)\neq 0\}$, where $\dom(f)$ denotes the domain of $f$. 
	We denote by $k\pow{G}$ the set of all total functions $s\colon G \to k$ whose support $\supp(s)$ is a well-ordered subset of $G$. For any $g\in G$, we denote by $t^g$ the characteristic function mapping $g$ to $1$ and everything else to $0$.
	This way we can express an element $s\in k\pow{G}$ as a power series $s=\sum_{g\in G}s_gt^g$, where $s_g=s(g)\in k$.
	For any power series $s,p\in k\pow{G}$, their sum is given by
{	\begin{equation}\notag\label{eq:sum-series}
		s+p=\sum_{g\in G}(s_g+p_g)t^g 
	\end{equation}}
	and their product by 
	{\begin{equation}\label{eq:prod-series}
		sp=\sum_{g\in G}q_gt^g  \text{ with } q_g=\sum_{h\in G}s_hp_{g-h}. 
	\end{equation}}
	Note that $ \supp(s+p)\subseteq \supp(s)\cup\supp(p) $ and $ \supp(sp)\subseteq\supp(s)+\supp(p) = \{g+h\mid g\in \supp(s), h\in \supp(p)\} $.
	These operations make $k\pow{G}$ a field (cf.~\cite{hahn,neumann}).
	We denote by $\cF_{\fin}$ the family of all finite subsets of $G$. 
	The \textbf{group ring} $ k[G] = \{s\in k\pow{G}\mid \supp(s)\in \cF_{\fin}\} $ is a subring of $ k\pow{G} $. Its fraction field $\ff(k[G])$ (within $k\pow{G}$) is denoted by $k(G)$.
	\begin{definition}\thlabel{def:hahn}
		We call $ k\pow{G} $ the \textbf{maximal Hahn field} (with coefficient field $k$ and value group $G$) and  
		$ k(G) $ the \textbf{minimal Hahn field}. Generally, we call a \textbf{Hahn field} any field $K$ such that $ k(G)\subseteq K \subseteq k\pow{G} $. 
	\end{definition} We define the \textbf{canonical valuation} on a Hahn field $K$ by
	$ v\colon K^\times\to G,\ s\mapsto \min\supp(s) $. This valuation has residue field $k$ and value group $G$.
	
	We use the following terminology for certain properties of $\cF$:
	\begin{itemize}
		\item $\cF$ is closed under subsets if for any $A\in \cF$ and any $B\subseteq A$ also $B\in \cF$.
		
		\item $\cF$ is closed under finite unions if for any $A,B\in \cF$ also $A\cup B \in \cF$.
		
		\item $\cF$ is closed under addition if for any $A,B\in \cF$ also $A+B=\{a+b\mid a\in A, b\in B\} \in \cF$.
		
		\item $\cF$ is closed under translation if for any $A\in \cF$ and any $g\in G$ also $A+g=\{a+g\mid a\in A\} \in \cF$.
	\end{itemize}	

\subsection{Kronecker's characterisation}

		Let $k\pow{t}=k\pow{\Z}$ denote the field of (formal) Laurent series over $k$ in one variable $t$. Note that with this notation we have $k(t)=k(\Z)$.
		We briefly present the main result of Kronecker~\cite[Section~VII\,ff.]{kronecker} exhibiting a way for determining whether a Laurent series in $k\pow{t}$ is already contained in  $k(t)$
		by means of linear recurrence relations (see also \cite[page~392\,f.]{hardy} or \cite[page~91\,ff.]{ribenboim} for further details in modern mathematical notation).
		
		\begin{fact}
			\thlabel{prop:laurentseries}
			Let $$s=\sum_{i\in \Z} s_it^i \in k\pow{t}.$$
			Then $s\in k(t)$ if and only if the following holds: there exist ${\ell}\in \N$ and $r_1,\ldots,r_{\ell}\in k$ such that for any $n> {\ell}$ the linear recurrence relation (on the coefficients of $s$)  with constant coefficients $r_1,\ldots,r_\ell$
			\begin{align}s_n=\sum_{j=1}^{{\ell}} r_j s_{n-j}\label{eq:linrec}\end{align}
			holds.
		\end{fact}

		\begin{remark}
			\begin{enumerate}[label = (\roman*)]
			\item In \Autoref{prop:laurentseries}, the element $s\in k(t)$ is not uniquely determined by $\ell$ and $r_1,\ldots,r_{\ell}$. For instance, we can consider $k=\Q$ and ${\ell}=r_1=1$. Then any $s={\sum_{i\leq 1}} s_it^i\in \Q(t)$, such as $s=1$ or $s=t^{-1}$,  satisfies (\autoref{eq:linrec}).
			
			Conversely, also ${\ell}$ and $r_1,\ldots,r_{\ell}$ are not uniquely determined by $s$. For instance, if we consider $s=t^{-1}\in \Q\pow{t}$, then we could also set ${\ell}=2$, $r_1=0$ and $r_2=1$. 
			
			\item {
			Given $\ell\in \N$, $r_1,\ldots,r_\ell\in k$, $i_0\in \Z$ with $i_0\leq 1$ and $s_{i_0},\ldots,s_{\ell}\in k$, there is a unique $s=\sum_{i=i_0}^\infty s_it^i\in k(t)$ satisfying (\autoref{eq:linrec}) for any $n>\ell$. The coefficients $s_{1},\ldots,s_{\ell}$ of $s$ may then be considered as initial values of $s$, and the following coefficients are determined by the linear recurrence relation. For instance, $s_{\ell+1}=r_1s_\ell+\ldots+r_\ell s_1$, so $s_{\ell+1}$ can be computed upon the input of constant coefficients $r_1,\ldots,r_\ell$ and initial values $s_{1},\ldots,s_{\ell}$.
			}
			\qed
			\end{enumerate}
		\end{remark}
		
		By setting $r_0=-1$, the linear recurrence relation (\autoref{eq:linrec}) can be reformulated as
		\begin{equation}
		0=\sum_{j=0}^\ell r_j s_{n-j}.
		\end{equation}
		We thus obtain the following consequence of \Autoref{prop:laurentseries}, which is a special case of the generalisation we will prove in \Autoref{lem:coefficientskG}.
		{For the convenience of the reader we include a  proof relying on \Autoref{prop:laurentseries}}.
		\begin{corollary}\thlabel{cor:laurentseries}
			Let $$s=\sum_{i\in \Z } s_it^i \in k\pow{t}.$$
			Then $s\in k(t)$ if and only if the following holds: there exist $\ell\in \N$, $g_0,\ldots,g_\ell\in \Z$ with $g_0<\ldots<g_\ell$ and $r_0,\ldots,r_\ell\in k$, not all equal to $0$, such that for any $n\in \Z\setminus\{g_0,\ldots,g_\ell\}$ the following linear recurrence relation with constant coefficients holds:
			\begin{align}
			\sum_{j=0}^\ell r_j s_{n-{g_j}}=0.\label{eq:linrec2}
			\end{align}
		\end{corollary}
	\begin{proof}
			{	
			Suppose first that $s=\sum_{i=i_0}^\infty s_it^i\in k(t)$, where $i_0\leq0$. Set $u_0=-1$. By \Autoref{prop:laurentseries}, 
			there exist ${m}\in \N$ and $u_1,\ldots,u_{m}\in k$ such that for any $n> {m}$:
			\begin{align}\sum_{j=0}^{{m}} u_j s_{n-j}=0\label{eq:newproof1}\end{align}
			Set $\ell=m-i_0$ and, for any $j\in \{0,\ldots,\ell\}$,  let $g_j=i_0+j$ and $$r_j=\begin{cases}0 & \text{if } g_j<0,\\
				u_{g_j} & \text{if } g_j\geq 0.
			\end{cases}$$
			Let $n\in \Z\setminus\{g_0,\ldots,g_\ell\}$. Then either $n>g_\ell=m$ or $n<g_0=i_0$. In the first case, due to (\ref{eq:newproof1}) we obtain
			\begin{align*}
				\sum_{j=0}^\ell r_j s_{n-g_j}=\sum_{j=-i_0}^\ell u_{i_0+j} s_{n-i_0-j}
				=\sum_{j=0}^{m} u_{j} s_{n-j}=0.
			\end{align*}
			In the second case, we have
			\begin{align*}
				\sum_{j=0}^\ell r_j s_{n-g_j}=\sum_{j=-i_0}^\ell u_{i_0+j} \underbrace{s_{n-i_0-j}}_{=0}
				 =0,
			\end{align*}
			as $n-i_0-j < i_0$ for $j\geq -i_0$.
			Thus, in either case (\ref{eq:linrec2}) holds.
			
			Conversely, suppose that there exist $\ell\in \N$, $g_0,\ldots,g_\ell\in \Z$ with $g_0<\ldots<g_\ell$ and $r_0,\ldots,r_\ell\in k$, not all equal to $0$, such that for any $n\in \Z\setminus\{g_0,\ldots,g_\ell\}$ (\ref{eq:linrec2}) holds. Let $j_0\in\{0,\ldots,\ell\}$ be least such that $r_{j_0}\neq 0$. Then for any $n>g_\ell$, we have
			\begin{align}\label{eq:newproof2}
				s_{n-g_{j_0}}=\sum_{j=j_0+1}^\ell \frac{-r_j}{r_{j_0}} s_{n-g_j}.
			\end{align}
			Set $m=|g_{j_0}|+|g_\ell|\in \N$ and let $n'>m$. Then $n'+g_{j_0}>g_\ell$. Applying (\ref{eq:newproof2}) yields
			\begin{align*}s_{n'}=s_{n'+g_{j_0}-g_{j_0}}=\sum_{j=j_0+1}^\ell \frac{-r_j}{r_{j_0}} s_{n'-(g_j-g_{j_0})}.\end{align*}
			As $0<g_{j_0+1}-g_{j_0}<\ldots<g_\ell-g_{j_0}$, we can thus choose suitable $u_1,\ldots,u_{m}\in k$ such that for any $n'> {m}$ 
			\begin{align*}s_{n'}=\sum_{j'=1}^{{m}} u_{j'} s_{n'-j'}\end{align*}
			holds. More precisely, for any $j'\in\{1,\ldots,m\}$ we set
			$$u_{j'}=\frac{-r_j}{r_{j_0}}$$ if 
			$j'=g_j-g_{j_0}$ for some $j\in\{j_0+1,\ldots,\ell\}$, and $u_{j'}=0$ otherwise. 
			By \Autoref{prop:laurentseries}, we obtain $s\in k(t)$, as required.
			}
	\end{proof}
		
		Note that the linear recurrence relation (\autoref{eq:linrec2}) in \Autoref{cor:laurentseries} needs to hold for all but finitely many coefficients, whereas the linear recurrence relation (\autoref{eq:linrec}) in \autoref{prop:laurentseries} only has to be satisfied for all coefficients whose index is greater than $l$. 
		The aforementioned results can be reformulated by means of Hankel matrices (see  \Autoref{fact:kroneckerhankel}).

\section{Generalised linear recurrence relations}\label{sec:gen-lin-rec}
	
	In this section, we introduce the main subject of study of this work: subsets of $k\pow{G}$ that are determined by generalised linear recurrence relations. We first set up the required concepts as well as notations and then establish several algebraic properties of sets determined by generalised linear recurrence relations depending on the underlying set of generalised linear recurrence sequences.
	
	\begin{definition}\thlabel{def:lrs}
		\begin{enumerate}[label = (\roman*)]
			\item A \textbf{(generalised) linear recurrence sequence} in $k\pow{G}$ is a partial function $r\colon G\rightharpoonup k$ whose domain $\dom(r)$ is a well-ordered subset of $G$. The order type of its domain corresponds to an ordinal number, which we denote by $\ot(r)$. The set of all linear recurrence sequences in $k\pow{G}$ is denoted by $\lc(k,G)$.
			
			\item Let $r\in \LRS(k,G)$. Then $r$ is called \textbf{$0$-free} if $r(g)\neq 0$ for any $g\in \dom(r)$ (i.e.\ $\supp(r)=\dom(r)$). We say that $r$ is \textbf{trivial} if $r(g)=0$ for any $g\in \dom(r)$.
			
			\item Let $r,r'\in \LRS(k,G)$. Then $r'$ is a \textbf{subsequence} of $r$ if $\dom(r')\subseteq \dom(r)$ and $r(g)=r'(g)$ for any $g\in \dom(r')$.
		\end{enumerate}	
	\end{definition}

	Note that any $r\in \LRS(k,G)$ can be decomposed into a trivial and a $0$-free subsequence as follows: Let $a,b\in \LRS(k,G)$ with $\dom(a)=\supp(r)$ and $\dom(b)=\dom(r)\setminus \supp(r)$. Moreover, set $a(g)=r(g)$ for any $g\in \dom(a)$ and $b(g)=0$ for any $g\in \dom(b)$. Then the graph of $r$ is the disjoint union of the graph of $a$ and the graph of $b$.

	\begin{notation}\thlabel{not:basicnotions}
		Let $r\in \lc(k,G)$ and let $s\in k\pow{G}$.
		\begin{enumerate}[label = (\roman*)]
			\item 	Let $\alpha=\ot(r)$. Then $r$ can be expressed as a sequence $(g_i,r_i)_{i<\alpha}$ in $G\times k$, where $(g_i)_{i<\alpha}$ is strictly increasing such that $\{g_i\mid i<\alpha\}=\dom(r)$ and $r_i=r(g_i)$ for any $i<\alpha$. 
			We then also write $r=(g_i,r_i)_{i<\alpha}$.

	\item We say that
			$$r^*=\sum_{g\in \dom(r)}r(g)t^g\in k\pow{G}$$ is the \textbf{power series associated} to $r$. Regarded as partial functions, $r^*$ is an extension of $r$ whose domain is $\dom(r^*)=G$. It has the property 
			$$r^*(g)=
			\begin{cases}
			r(g)&\text{ for }g\in \dom(r),\\
				0&\text{ for }g\in G\setminus \dom(r).
			\end{cases}$$
			Thus, $\supp(r^*)=\dom(r)\setminus r^{-1}(0)$. Moreover, $\supp(r^*)=\supp(r)$.
	\item We say that $ s^* \in \LRS(k,G) $ defined by
	$$
	s^*\colon \supp(s) \to k, 
	g\mapsto 
		s(g)$$ is the \textbf{linear recurrence sequence associated} to $s$.
		Note that $s^*$ is $0$-free and $\supp(s^*)=\supp(s)$.

	\item We write $r^{**}$ for $(r^*)^*$ and $s^{**}$ for $(s^*)^*$. Note that $s^{**}=s$. Note further that $r^{**}$ is the restriction of $r$ to the subset $\supp(r)$ of $\dom(r)$. Thus, $r^{**}$ is the largest {$0$-free subsequence} of $r$.
	\end{enumerate}
	\end{notation}

	\begin{definition}\thlabel{def:determined}
		Let $r=(g_i,r_i)_{i<\alpha}\in \lc(k,G)$. We define $\gen{r}\subseteq k\pow{G}$ to be the set of all power series $s\in k\pow{G}$  such that for any
		$$h\in G\setminus\dom(r)$$
		the following \textbf{(generalised) linear recurrence relation} holds:
		\begin{align}\sum_{i<\alpha}r_is_{h-g_i}=0.\label{eq:genlinrecrel}\end{align}
		(Note that $h-g_i$ is strictly decreasing, whence the sum in (\autoref{eq:genlinrecrel}) is finite.
		)
		We say that $r$ \textbf{determines} $\gen{r}$. For a set $R\subseteq \LRS(k,G)$ we set $$\gen{R}=\bigcup_{r\in R}\gen{r}$$
		and say that $R$ determines $\gen{R}$.
		For a set $S\subseteq k\pow{G}$ we say that $S$ is \textbf{determined} by (generalised) linear recurrence relations if there exists $R\subseteq \LRS(k,G)$ such that $\gen{R}=S$.
	\end{definition}

	\begin{example}
		Let $r=(g_i,r_i)_{i<\alpha}\in \lc(k,G)$. If $r$ is trivial, then (\autoref{eq:genlinrecrel}) always becomes $0=0$, i.e.\ the linear recurrence relation holds for any $s\in k\pow{G}$. Hence, $\gen{r}=k\pow{G}$, i.e.\ $r$ determines the maximal Hahn field $k\pow{G}$. 
		Conversely, if $r$ is non-trivial, then $r_j\neq 0$ for some $j<\alpha$. Let $h\in G$ with $h<g_0$ (and thus $h\in G\setminus \dom(r)$) and let $s\in k\pow{G}$ be the monomial $s=r_j^{-1}t^{h-g_j}$.
		Then $$\sum_{i<\alpha}r_is_{h-g_i}=r_js_{h-g_j}=r_jr_j^{-1}=1.$$ Hence, $s$ does not satisfy (\autoref{eq:genlinrecrel}), yielding $s\notin \gen{r}$.
		This shows that $\gen{r}\subsetneq k\pow{G}$, i.e.\ $r$ does not determine $k\pow{G}$.

		From the above, we obtain that for any $R\subseteq \lc(k,G)$ containing a trivial linear recurrence sequence, we have $\gen{R}=k\pow{G}$. \qed
	\end{example}

	For any non-trivial linear recurrence sequence $r$, its associated power-series $r^*$ is non-zero. In this case, the following lemma gives a description of the set determined by $r$.

	\begin{lemma}[Main Lemma]\thlabel{prop:description}
		Let $r\in \LRS(k,G)$ be non-trivial. Then
		\begin{align}\gen{r}=\setbr{\tfrac{b}{r^*} \mid  b\in k\pow{G}\!, \supp(b)\subseteq \dom(r)}.\label{eq:description}\end{align}
	\end{lemma}

	\begin{proof}
		Set $r=(g_i,r_i)_{i<\alpha}$.
		Let $b \in k\pow{G}$ with $\supp(b)\subseteq \dom(r)$ and set $a= \tfrac{b}{r^*}$. Further, let $h\in G\setminus\dom(r)$. Then {by the product formula \eqref{eq:prod-series}}
		$$\sum_{i<\alpha}r(g_i)a_{h-g_i}=(r^*a)_h=b_h=0,$$
		as $h\notin \dom(r)\supseteq \supp(b)$. Hence, $a\in \gen{r}$, showing that $$\gen{r}\supseteq\setbr{\tfrac{b}{r^*} \mid  b\in k\pow{G}\!, \supp(b)\subseteq \dom(r)}.$$
		For the converse, let $a\in k\pow{G}$ with $a\in \gen{r}$. Set
		$$b_h:=(r^*a)_h=
		\sum_{i<\alpha}r(g_i)a_{h-g_i}$$
		for any $h\in G$.
		Then for any $h\in G\setminus \dom(r)$ we have $b_h=0$, as $a\in \gen{r}$. Hence, $\supp(b)\subseteq \dom(r)$. Since for any $h \in G$ we have $b_h=(r^*a)_h$, we obtain $b=r^*a$. Hence, $a=\tfrac b{r^*}$, as required.
	\end{proof}

		Since (\autoref{eq:description}) is a useful and important characterisation of sets determined by a single non-trivial linear recurrence sequence, we will use it throughout this work without explicitly referring to \Autoref{prop:description} each time.

	\begin{example}
		Consider
				$$s = 1 + 1t^{\sqrt{2}} + 2t^{2\sqrt{2}} + 3t^{3\sqrt{2}} + 5t^{4\sqrt{2}} + 8t^{5\sqrt{2}} + 13t^{6\sqrt{2}}+ \ldots \in k\pow{\R}.$$
			For any $h\in \R\setminus\{0,\sqrt{2},2\sqrt{2}\}$, the linear recurrence relation
			$$s_h-s_{h-\sqrt{2}}-s_{h-2\sqrt{2}}=0$$
			holds. 
			Thus, $s\in \setbr{\left.\!\frac{a+bt^{\sqrt{2}}+ct^
					{2\sqrt{2}}}{1-t^{\sqrt{2}}-t^{2\sqrt{2}}}\ \right|\ a,b,c\in k}$. 
			Indeed,
			$$s=\frac{1}{1-t^{\sqrt{2}}-t^{2\sqrt{2}}}\in k(\R).$$\qed
	\end{example}

	\begin{proposition}\thlabel{cor:description}
		Let $r\in \LRS(k,G)$. Then $\gen{r}\subseteq k\pow{G}$ is a vector space over $k$ containing $k$.
	\end{proposition}

	\begin{proof}
		If $r$ is trivial, then $\gen{r}=k\pow{G}$, which is even a field extension of $k$. Otherwise, for any $a,b\in k\pow{G}$ with $\supp(a),\supp(b)\subseteq\dom(r)$ and any $\lambda \in k$ we have
		$$\frac{a}{r^*}+\lambda\cdot  \frac{b}{r^*}=\frac{a+\lambda b}{r^*}.$$
		Now $\supp(a+\lambda b)\subseteq \supp(a)\cup\supp(b)\subseteq \dom(r)$. Hence, also $\frac{a+\lambda b}{r^*}$ lies in $\gen{r}$.
		
		Finally, note that $\supp(r^*)\subseteq \dom(r)$, as $\supp(r^*)=\dom(r)\setminus r^{-1}(0)$. Hence, $1=\frac{r^*}{r^*}\in \gen{r}$. Since $\gen{r}$ is closed under multiplication with elements from $k$, we obtain $k\subseteq \gen{r}$.
	\end{proof}
	
	While \Autoref{cor:description} shows that a single linear recurrence sequence determines a $k$-vector subspace of $k\pow{G}$, a non-empty set $R$ of linear recurrence sequences does not even necessarily determine a subset of $k\pow{G}$ that is closed under addition, as the following example shows.
	
	\begin{example}\thlabel{ex:additionclosed}
		Consider the power series $a=1+t$ and $b=1-t$ in $\Q\pow{\Z}$ and let $R=\{a^*,b^*\}$ be the set of associated linear recurrence sequences. Then
		\begin{align*}
			\gen{R}=&\{\tfrac ca\mid c\in \Q\pow{\Z}, \supp(c)\subseteq\{0,1\}\}\cup\{\tfrac db \mid d\in \Q\pow{\Z}, \supp(d)\subseteq\{0,1\}\}\\
			&=\{\tfrac{c_0+c_1t}{1+t} \mid c_0,c_1\in \Q\}\cup\{\tfrac{d_0+d_1t}{1-t} \mid d_0,d_1\in \Q\}.
		\end{align*}
		Thus, we have $\tfrac{1}{1+t},\tfrac{1}{1-t}\in \gen{R}$. However, $$\tfrac{1}{1+t}+\tfrac{1}{1-t}=\tfrac{2}{1-t^2}\notin\gen{R},$$
		as neither the equation $\tfrac{2}{1-t^2}=\tfrac{c_0+c_1t}{1+t}$ nor the equation $\tfrac{2}{1-t^2}=\tfrac{d_0+d_1t}{1-t}$ has a solution over $\Q$. Hence, $\gen{R}$ is not closed under addition.
		
		Note that for a similar reason we have $\tfrac{1}{1+t}\cdot\tfrac{1}{1-t}=\frac{1}{1-t^2}\notin \gen{R}$, whence  $\gen{R}$ is also not closed under multiplication.  \qed
	\end{example}

	We now establish sufficient conditions for a set determined by linear recurrence relations to be closed under addition and multiplication as well as taking multiplicative inverses of non-zero elements.

	\begin{lemma}\thlabel{lem:closure1}
		Let $R\subseteq \LRS(k,G)$. Suppose that for any $r_1,r_2\in R$ the following $r\in \LRS(k,G)$ with domain\footnote{Note that $\dom(r)$ is well-ordered, as the sum of two well-ordered subsets of $G$ is also well-ordered.} $\dom(r)=\dom(r_1)+\dom(r_2)$ is also {an element of} $R$: for any $g\in \dom(r_1)+\dom(r_2)$ set
		$$r(g):=\begin{cases}
		(r_1^*\cdot r_2^*)(g)& \text{ if } g\in \supp(r_1^*r_2^*),\\
		0 & \text{ if } g\in (\dom(r_1)+\dom(r_2))\setminus \supp(r_1^*r_2^*).
		\end{cases}$$	
		Then $\gen{R}$ is closed under addition and multiplication.
	\end{lemma}

	\begin{proof}
		Let $r_1,r_2\in R$, let $a_1\in \gen{r_1}$ and let $a_2\in \gen{r_2}$. 
		We may assume that $r_1$ and $r_2$ are non-trivial, as otherwise we already have $\gen{R}=k\pow{G}$. 
		Fix $b_1,b_2\in k\pow{G}$ with $a_1=\tfrac{b_1}{r_1^*}$ and $a_2=\tfrac{b_2}{r_2^*}$, where $\supp(b_1)\subseteq \dom(r_1)$ and $\supp(b_2)\subseteq \dom(r_2)$. Then
		$$a_1+a_2=\tfrac{b_1r_2^*+b_2r_1^*}{r_1^*r_2^*}\text{\ \ and\ \ }a_1a_2 = \tfrac{b_1b_2}{r_1^*r_2^*}.$$
		Note that by definition of $r$, we have $r^*=r_1^*r_2^*$, as $\supp(r_1^*r_2^*)\subseteq \supp(r_1^*)+\supp(r_2^*)\subseteq \dom(r_1)+\dom(r_2)$.
		Thus, 		
		$$a_1+a_2=\tfrac{b_1r_2^*+b_2r_1^*}{r^*}\text{\ \ and\ \ }a_1a_2 = \tfrac{b_1b_2}{r^*}.$$
		Moreover,
		\begin{align*}
			\supp(b_1r_2^*+b_2r_1^*)&\subseteq \supp(b_1r_2^*)\cup \supp(b_2r_1^*)\\
			&\subseteq [\supp(b_1)+\supp(r_2^*)]\cup [\supp(b_2)+\supp(r_1^*)]\\
			&\subseteq [\dom(r_1)+\dom(r_2)]\cup [\dom(r_2)+\dom(r_1)]\\
			&= \dom(r)
		\end{align*}
		and
		\begin{align*}
		\supp(b_1b_2)\subseteq  \supp(b_1)+\supp(b_2)\subseteq\dom(r_1)+\dom(r_2)=\dom(r).
		\end{align*}
		Hence, $$a_1+a_2,a_1a_2\in \gen{r}\subseteq\gen{R}.$$
	\end{proof}

	\begin{lemma}\thlabel{lem:closure2}
		Let $R\subseteq\LRS(k,G)$. Suppose that for any $r\in R$ any non-trivial $a\in \LRS(k,G)$ with domain $\dom(a)=\dom(r)$ is also {an element of} $R$. Then $\gen{R}$ is closed under taking multiplicative inverses of non-zero elements.
	\end{lemma}

	\begin{proof}
		Let $r\in R$. If $r$ is trivial, then $\gen{R}=k\pow{G}$ and we are done. Otherwise, let $s\in \gen{r}\setminus\{0\}$ and let $b\in k\pow{G}$ with $\supp(b)\subseteq \dom(r)$ and
		$$s=\tfrac{b}{r^*}.$$
		Let $a\in R$ be the linear recurrence sequence with domain $\dom(r)$ given by
		$$a(g):=\begin{cases}
		b(g)&\text{ for }g\in \supp(b),\\
		0&\text{ for }g\in \dom(r)\setminus\supp(b).
		\end{cases}$$
		Then $a^*=b$. Since $\supp(r^*)\subseteq \dom(r)=\dom(a)$, we obtain
		$$s^{-1}=\tfrac{r^*}{b}=\tfrac{r^*}{a^*}\in \gen{a}\subseteq \gen{R},$$
		as required.
	\end{proof}

	Combining \Autoref{lem:closure1} and \Autoref{lem:closure2}, we obtain the following sufficient conditions on a set $R\subseteq\LRS(k,G)$ in order that $\gen{R}$ is a subfield of $k\pow{G}$.

	\begin{theorem}\thlabel{prop:field}
		Let $R\subseteq \LRS(k,G)$ be non-empty. Suppose that for any $r_1,r_2\in R$ every non-trivial $a,b\in \LRS(k,G)$ with $\dom(a)=\dom(r_1)$ and $\dom(b)=\dom(r_1) + \dom(r_2)$ is also {an element of} $R$. Then $\gen{R}$ is a subfield of $k\pow{G}$ containing $k$.
	\end{theorem}

	\begin{proof}
		Since $R$ is non-empty, \Autoref{cor:description} shows that $k\subseteq \gen{R}$.\linebreak
		\Autoref{lem:closure2} immediately implies that $\gen{R}$ is closed under taking multiplicative inverses of non-zero elements. 
		If $R$ contains a trivial sequence, then $\gen{R}=k\pow{G}$. Otherwise, the sequence $r$ with $\dom(r)=\dom(r_1)+\dom(r_2)$ in	\Autoref{lem:closure1} is non-trivial for any $r_1,r_2\in R$ and thus contained in $R$. This implies that $\gen{R}$ is closed under addition and multiplication.
	\end{proof}

\section{Linear recurrence Hahn fields}\label{sec:hahnfields}

	\Autoref{prop:field} enables us to consider specific subfields of $k\pow{G}$ determined by linear recurrence relations. More specifically, if $R\subseteq \LRS(k,G)$ consists of \emph{all} non-trivial linear recurrence sequences with a particular domain, then $\gen{R}$ forms a subfield of $k\pow{G}$. We make this precise in the following. Recall that $\cF$ always denotes a family of well-ordered subsets of $G$.
	
	\begin{definition}\thlabel{def:fseq}
		We denote by $S(\cF)$ the set of all non-trivial $r\in \LRS(k,G)$ with $\dom(r)\in\cF$. A linear recurrence sequence $r\in S(\cF)$ is called an \textbf{$\cF$-sequence}.
	\end{definition}

	For later use, we observe in the following proposition that if $\cF$ contains an upper bound $A\in \cF$ with respect to its partial ordering $\subseteq$, then any power series determined by an $\cF$-sequence is already determined by an $\{A\}$-sequence.
		
	\begin{proposition}\thlabel{prop:upperbound}
		Suppose that there is some $A\in \cF$ such that $\cF\subseteq \mathcal{P}(A)$. Then $\gen{S(\cF)}=\gen{S(\{A\})}$.
	\end{proposition}

	\begin{proof}
		Since $\{A\}\subseteq \cF$, we have $\gen{S(\{A\})}\subseteq \gen{S(\cF)}$. For the converse, let $r$ be an $\cF$-sequence and let $b\in k\pow{G}$ with $\supp(b)\subseteq \dom(r)\in \cF$. We have to show that $\tfrac b {r^*}\in \gen{S(\{A\})}$. Since $\dom(r)\in \cF$, we have $\dom(r)\subseteq A$. Let $a$ be the $\{A\}$-sequence given by $a(g)=r(g)$ for $g\in \dom(r)$ and $a(g)=0$ for $g\in A\setminus \dom(r)$. Then $a^*=r^*$ and $\supp(b)\subseteq A=\dom(a)$, showing that $\tfrac b {r^*}=\tfrac b {a^*}\in \gen{S(\{A\})}$.
	\end{proof}
		
	We now turn to fields determined by $\cF$-sequences. 
	\Autoref{prop:field} applied to $R=S(\cF)$ immediately implies the following.
	
	\begin{corollary}\thlabel{cor:field}
		Suppose that $\cF\notin \{\emptyset,\{\emptyset\}\}$ is closed under addition. Then $\gen{S(\cF)}$ is a subfield of $k\pow{G}$ containing $k$.
	\end{corollary}

	Recall that $\cF_{\fin}$ denotes the set of all finite subsets of $G$. Hence, the set of $\cF_{\fin}$-sequences $S(\cF_{\fin})$ consists of all non-trivial linear recurrence sequences with finite domain.
	
	\begin{proposition}\thlabel{prop:finitesets}
		The minimal Hahn field $k(G)$ is determined by $S(\cF_{\fin})$, i.e.\ $\gen{S(\cF_{\fin})}=k(G)$.
	\end{proposition}
	
	\begin{proof}
		By \Autoref{cor:field}, $\gen{S(\cF_{\fin})}\subseteq k\pow{G}$ is a field containing $k$, completing the proof for the case $G=\{0\}$. Otherwise, let $g\in G\setminus \{0\}$ and let $r$ be the $\cF_{\fin}$-sequence with $\dom(r)=\{0,g\}$, $r(0)=1$ and $r(g)=0$. Then $$t^g=\tfrac{t^g}{1}=\tfrac{t^g}{r^*}\in \gen{r}\subseteq \gen{S(\cF_{\fin})}.$$
		This shows that $k(G)\subseteq \gen{S(\cF_{\fin})}$.
		
		Conversely, let $s\in \gen{S(\cF_{\fin})}$. Then there exist $r\in S(\cF_{\fin})$ and  $b\in k\pow{G}$ with $\supp(b)\subseteq \dom(r)$ and  $s=\tfrac b{r^*}$. Since $r$ has a finite domain, both $r^*$ and $b$ have finite support. Hence, $s\in k(G)$, as required. 
	\end{proof}
	
	\Autoref{cor:field} and \Autoref{prop:finitesets} imply the following.
	
	\begin{corollary}\thlabel{cor:hahnfield}
		Suppose that $\cF_{\fin}\subseteq \cF$. 
		Then $k(G)\subseteq \gen{S(\cF)}$.
		In particular, 
		\begin{itemize}
			\item if $\gen{S(\cF)}$ is a field, then $ \gen{S(\cF)} $ is a Hahn field, and
			\item if $ \cF $ is closed under addition, then $ \gen{S(\cF)} $ is a Hahn field.
		\end{itemize}
	\end{corollary}
	
	\Autoref{cor:hahnfield} can be used to construct an example of a Hahn field that is not determined by linear recurrence relations as follows.
	
	\begin{lemma}\thlabel{prop:uncountable}
		Let $R\subseteq \LRS(k,G)$ such that $k(G)\subsetneq \gen{R}$. Then $\gen{R}$ is uncountable. In particular, any non-minimal Hahn field that is determined by linear recurrence relations is uncountable.
	\end{lemma}
	
	\begin{proof}
		First note that $R$ is non-empty and thus $G\neq \{0\}$, as otherwise $k(G)=k=\gen{R}$ by \Autoref{cor:description}. 
		If $R$ contains a trivial linear recurrence sequence, then $\gen{R}=k\pow{G}$, which is always uncountable. Otherwise, 
		\Autoref{prop:finitesets} implies that $R$ contains a non-trivial $r\in \LRS(k,G)$ with infinite domain. Now $$\{\tfrac{a}{r^*}\mid a\in k\pow{G}\!, \supp(a)\subseteq \dom(r)\}=\gen{r}\subseteq \gen{R}.$$
		Hence, for every $B\subseteq \dom(r)$, the set $\gen{r}$ contains the power series
		$$\tfrac{\sum_{g\in B}t^g}{r^*}.$$
		Thus, $\gen{r}$ contains at least $|\cP(\dom(r))|=2^{|\dom(r)|}\geq 2^{\aleph_0}$ elements, whence $\gen{R}$ is uncountable.
	\end{proof}
	
	While \Autoref{prop:finitesets} shows that any minimal Hahn field is determined by linear recurrence relations, and thus also every countable minimal Hahn field such as $\Q(\Z)$, any countable Hahn field that is \emph{not} the minimal Hahn field is \emph{not} determined by linear recurrence relations, as \Autoref{prop:uncountable} shows. For instance, let $C$ be the relative algebraic closure of $\Q(\Z)$ inside $\Q\pow{\Z}$. Since $\sqrt{1+t}\in C\setminus \Q(\Z)$, we have that $C$ is a countable non-minimal Hahn field
	and thus not determined by linear recurrence relations.
	
	We now come to the generalisation of \Autoref{cor:laurentseries} to fields of generalised power series.
	
	\begin{theorem}\thlabel{lem:coefficientskG} 
		Let $s=\sum_{g\in G}s_gt^g\in k\pow{G}$. Then $s\in k(G)$ if and only if the following holds:
		there exist $\ell\in \N$, finite sets $A=\{g_0,\ldots,g_\ell\}\subseteq G$ with $g_0<\ldots<g_\ell$ and $\{r_0,\ldots,r_\ell\}\subseteq k$, not all $r_i$ equal to $0$, such that for any $h\in G\setminus A$ the following linear recurrence relation holds:
		\begin{align}\sum_{j=0}^\ell r_js_{h-g_j}=0\label{eq:genrecrel}\notag\end{align}
	\end{theorem}
	
	\begin{proof} 
		This is basically a consequence of \Autoref{prop:finitesets}: 
		Suppose that  $s\in k(G)=\gen{S(\cF_{\fin})}$. Then there exist $\ell\in\N$, $A= \{g_0,\ldots,g_\ell\}\in \cF_{\fin}$ and a non-trivial linear recurrence sequence $r=(g_i,r_i)_{i<\ell+1}$ in $S(\cF_{\fin})$ with $s\in \gen{r}$. Hence,  for any $h\in G\setminus \dom(s)=G\setminus A$ we have
		\begin{align*}\sum_{j=0}^\ell r_js_{h-g_j}=0.\end{align*}
		Conversely, suppose that there exist $\ell\in\N$, finite sets $A=\{g_0,\ldots,g_\ell\}\subseteq G$ with $g_0<\ldots<g_\ell$ and $\{r_0,\ldots,r_\ell\}\subseteq k$, not all $r_i$ equal to $0$, such that for any $h\in G\setminus A$ we have
		\begin{align*}\sum_{j=0}^\ell r_js_{h-g_j}=0.\end{align*}
		Then $s\in \gen{r}\subseteq \gen{S(\cF_{\fin})}=k(G)$ for $r=(g_i,r_i)_{i<\ell+1}\in S(\cF_{\fin})$.
	\end{proof}

\section{Generalised Hankel matrices}
\label{sec:gen-Hank-Matrices}

	We now generalise Kronecker's criterion on the rationality of power series $s\in k\pow{t}$ in terms of finite ranks of their Hankel matrix (see \cite[Section~IX]{kronecker}) to generalised power series $s\in k\pow{G}$ for an arbitrary ordered abelian group $G$. {\Autoref{fact:kroneckerhankelgen}, the main result of this section, applies to power series of the form that is explicitly given in \Autoref{rmk:hankelform}.} Generally, a Hankel matrix is a possibly infinite matrix with constant {ascending diagonals}. In the following, we make precise how to obtain a Hankel matrix associated to a given power series $s \in k\pow{G}$ with order type at most $\omega$ and fixed increment $g\in G$ {(see \Autoref{def:Hankelclassic}~\ref{def:Hankelclassic:2})}.
	{ In our context Hankel matrices are always of dimension $\omega\times\omega$, also if $\supp(s)$ is finite. }
	
	\begin{definition}\thlabel{def:Hankelclassic}
		\begin{enumerate}[label = (\roman*)]
			\item Let $g_0,g\in G$ with $g>0$ and let $r = (g_0+ig,r_i)_{i<\omega}\in \LRS(k,G)$.
			Then the \textbf{Hankel matrix} associated to $ r $ is given by
			$\Hk{r}=(\Hk{r}_{ij})_{i,j\in\N} \in k^{\N\times \N} $ with 
			$$ \Hk{r}_{ij}:= r_{i+j}.$$ for any $i,j\in \N$.
			
			\item\label{def:Hankelclassic:2}{ Let $s\in k\pow{G}\setminus\{0\}$ satisfy the following condition: there exists $g\in G^{>0}$ such that for $r_{s,g} := (g_i,s_{g_i})_{i<\omega}\in \LRS(k,G)$ with $(g_i)_{i<\omega}:=(v(s)+ig)_{i<\omega}$ we have $r_{s,g}^*=s$.
			Then the \textbf{Hankel matrix} associated to $ s $ with \textbf{increment} $g$ is given by
			$\Hk{s,g}:=\Hk{r_{s,g}}$.}
		\end{enumerate}		
	\end{definition}
	
	\begin{remark}\label{rmk:hankelform}
		Let $s\in k\pow{G}\setminus\{0\}$, $g\in G^{>0}$ and let $r_{s,g} = (g_i,s_{g_i})_{i<\omega}\in \LRS(k,G)$ be as in \Autoref{def:Hankelclassic}~\autoref{def:Hankelclassic:2}. Since $$r_{s,g}^*=\sum_{i<\omega}s_{g_i}t^{g_i},$$ we have that $r_{s,g}^*=s$ if and only if $\dom(r_{s,g})=\{v(s)+ig\mid i<\omega\}$ contains $\supp(s)$. Hence,  the Hankel matrix associated to $ s $ with increment $g$ exists if and only if $s$ is of the form
		$$s=t^{v(s)}\sum_{i<\omega}s_{v(s)+ig}t^{ig}.$$
		{ 
		Setting $h=v(s)$ and $F(x)=\sum_{i=0}^\infty s_{h+ig}x^i\in k\rpow{x}$, we can express such a power series $s$ as
		$$s=t^hF(t^g).$$
		}
		\qed
	\end{remark}

	\begin{example}\thlabel{ex:hankelexp}
		Let $s\in k\pow{t}\setminus\{0\}$. Then $r_{s,1} = (v(s)+i,s_{g_i})_{i<\omega}\in \LRS(k,G)$ satisfies $r_{s,1}^*=s$. The Hankel matrix $\Hk{s,1}=\Hk{r_{s,1}}$ associated to $s$ with increment $1$ is the original Hankel matrix that Kronecker associated to $s$.
		
		For instance, let $k$ be of characteristic $0$ and consider 
		$$e^t:=\sum_{i=0}^\infty\frac{t^i}{i!} \in k\pow{t}.$$
			Then $r_{e^t,1}=(i,\tfrac{1}{i!})_{i<\omega}\in \LRS(k,\Z)$ and 
			\[	\Hk{e^t,1} = \begin{pmatrix}
			\frac{1}{0!} & \frac{1}{1!} & \frac{1}{2!} &\reflectbox{$\ddots$}\\[2pt]
			\frac{1}{1!} & \frac{1}{2!} & \frac{1}{3!} &\reflectbox{$\ddots$}\\[2pt]
			\frac{1}{2!} & \frac{1}{3!} & \frac{1}{4!}&\reflectbox{$\ddots$} \\
			\reflectbox{$\ddots$} &\reflectbox{$\ddots$} &\reflectbox{$\ddots$} &
			\end{pmatrix}\in k^{\N\times\N}.
			\]
		If we choose the increment to be $g=\tfrac 12$, then the Hankel matrix of $e^t$ becomes
		\[	\Hk{e^t,\frac 12} = \begin{pmatrix}
		\frac{1}{0!} & 0 & \frac{1}{1!} &\reflectbox{$\ddots$}\\[2pt]
		0 & \frac{1}{1!} & 0 &\reflectbox{$\ddots$}\\[2pt]
		\frac{1}{1!} & 0 & \frac{1}{2!}&\reflectbox{$\ddots$} \\
		\reflectbox{$\ddots$} &\reflectbox{$\ddots$} &\reflectbox{$\ddots$} &
		\end{pmatrix}\in k^{\N\times\N}.
		\]
	\qed
	\end{example}

	The following is Kronecker's original result (see also \cite[Section~I.3, Lemma III]{salem} for a proof).

	\begin{fact}\thlabel{fact:kroneckerhankel}
		Let $ s\in k\pow{t}\setminus\{0\} $. Then 
		$ s\in k(t) $ if and only if the rank of $ \Hk{s,1} $ is finite. 
	\end{fact}

	We now generalise \Autoref{fact:kroneckerhankel} to Hahn fields with any value group $G$.

	\begin{theorem}\thlabel{fact:kroneckerhankelgen}
		Let $s\in k\pow{G}\setminus\{0\}$, let $g\in G^{>0}$ and suppose that $r_{s,g}^*=s$.
		Then 
		$ s\in k(G) $ if and only if the rank of the Hankel
		matrix $ \Hk{s,g} $  is finite.
	\end{theorem}
	
	\begin{proof}
		Since $G$ is an ordered abelian group and $g>0$, we can fix an embedding $\Z\hookrightarrow G, 1\mapsto g$ of ordered groups and identify $\Z$ with its image in $G$ via this map, i.e.\ we regard $ng$ as $n\in \Z$ for any $n\in \Z$. Moreover, we obtain an embedding of fields $k\pow{t}=k\pow{\Z}\subseteq k\pow{G}$ via this identification.
		Now $\Hk{s,g}=\Hk{s,1}$, and \Autoref{fact:kroneckerhankel} implies that the rank of $\Hk{s,g}$ is finite if and only if $s\in k(t)$. The latter holds if and only if $s=\tfrac{p(t)}{q(t)}$ for some $p,q\in k[t]$, which, by our identification, is equivalent to $s=\tfrac{p(t^{g})}{q(t^{g})}\in k(G)$.
	\end{proof}

	As an application of \Autoref{fact:kroneckerhankelgen} we show in the following that for $g>0$ the power series $e^{t^g}=\sum_{n=0}^\infty\frac{t^{ng}}{n!}$ is never contained in the minimal Hahn field $k(G)$.
			
	\begin{lemma}\thlabel{lem:ax}
		Let $ \Char(k) = 0 $. Then
		$ e^{t} $ is transcendental over $ k(t) $.
	\end{lemma}
	\begin{proof}
		This is a special case of \cite[Corollary~1]{ax}, where we set $ C = k $, $ r = n = 1 $ and $ y_1 = t_1 = t $.
	\end{proof}
		
	\begin{example}\label{ex:k(G)-not-k-hull0}
		Suppose that $\Char(k)=0$ and consider the power series 
		$$  E=e^{t^{g}}:=\sum_{{i}=0}^\infty\frac{t^{{i}g}}{{i}!}.  $$ We show that $E\notin k(G)$.
		Note that $\Hk{E,g}=\Hk{e^{t},1}$ (see \Autoref{ex:hankelexp}). Since by \Autoref{lem:ax}, we have that $e^t\notin k(t)$, \Autoref{fact:kroneckerhankel} shows that the rank of $\Hk{e^{t},1}$ is infinite. Hence, also the rank of $\Hk{E,g}$ is infinite, and $E\notin k(G)$ by \Autoref{fact:kroneckerhankelgen}. \qed
	\end{example}

\section{Relation to Rayner fields}\label{sec:rayner}
	
	Recall that $\cF$ denotes a family of well-ordered subsets of $G$.
	Following the terminology of \cite[Definition~2.4]{kks}, we call $$k\pow{\cF}=\{a\in k\pow{G}\mid \supp (a) \in \cF\}$$ the \textbf{$k$-hull}\label{page:khull} of $\cF$ in $k\pow{G}$.
	If $\cF$ consists of all well-ordered subsets of some set $A\subseteq G$, then we simply denote $k\pow{\cF}$ by $k\pow{A}$.
	
	In \cite{kks}, we studied in detail how conditions on $\cF$ relate to set theoretic and algebraic properties of its $k$-hull. A particular focus was laid on Rayner fields: a $k$-hull $k\pow{\cF}$ is a \textbf{Rayner field} if $\cF$ is non-empty, is closed under subsets, finite unions and translations, $\Span_\Z\!\brackets{\bigcup_{A\in \cF}A}=G$ and for any $A\in \cF$ with $A\subseteq G^{\geq 0}$ also $\Span_{\N}(A)\in \cF$ (see \cite[Definition~3.2]{kks}).
	We showed in \cite[Lemma~3.17]{kks} that any Rayner field is a Hahn field and in \cite[Theorem~3.18]{kks} that in the case $\Char(k)=0$ any $k$-hull that is a Hahn field is also a Rayner field. This section explores which Rayner fields and, more generally, $k$-hulls are determined by linear recurrence relations.

	We first observe that not all $k$-hulls  are determined by linear recurrence relations: Choose $\cF\neq \emptyset$ such that it does not contain $\{0\}$. Then $k\not\subseteq k\pow{\cF}$, whence $k\pow{\cF}$ is not determined by linear recurrence relations by \Autoref{cor:description}. However, under certain conditions on $\cF$, we can always find a set of linear recurrence sequences $R_{\cF}$ such that $\gen{R_{\cF}}=k\pow{\cF}$.

	\begin{notation}\thlabel{not:rayner}
			Let $A\subseteq G$ be well-ordered. We denote by $\cF_A$ the family consisting of all subsets of $A\cup\{0\}$, i.e.\ $\cF_A=\mathcal{P}(A\cup\{0\})$. Moreover, we let $r_A\in \LRS(k,G)$ with domain $\dom(r_A)=A\cup\{0\}$ mapping $0$ to $1$ and everything else to $0$, and set $R_\cF=\{r_A\mid A\in \cF\}$.
	\end{notation}

	Note that for any $A\in \cF$, we have $r_A^*=1\in k\pow{G}$. 

	\begin{lemma}\thlabel{lem:rayner1}
		Let $A\subseteq G$ be well-ordered. Then 
		$\gen{r_A}=k\pow{\cF_A}$.
	\end{lemma}

	\begin{proof}
		Let $a\in \gen{r_A}$. Then there exists $b\in k\pow{G}$ with $\supp(b)\subseteq \dom(r_A)=A\cup\{0\}$ such that $a=\tfrac{b}{r_A^*}=b$. Hence, $a=b\in k\pow{\cF_A}$.
		
		Conversely, let $a\in k\pow{\cF_A}$. Then $\supp(a)\in \cF_A$ and thus $\supp(a)\subseteq A\cup \{0\}=\dom(r_A)$. Hence, $a=\tfrac{a}{r_A^*}\in \gen{r_A}$.
	\end{proof}

	\begin{proposition}\label{prop:rayner2}
		Suppose that $\cF$ is closed under subsets and unions with $\{0\}$ (i.e.\ $A\in \cF$ implies $A\cup\{0\}\in \cF$). Then 
		$$\gen{R_{\cF}}=k\pow{\cF}.$$
	\end{proposition}

	\begin{proof}
		We may suppose that $\cF$ is non-empty, as otherwise $\gen{R_{\cF}}=\emptyset =k\pow{\cF}$. 
		Let $A\in \cF$. Then $\cF_A\subseteq \cF$. \Autoref{lem:rayner1} implies $\gen{r_A}=k\pow{\cF_A}\subseteq k\pow{\cF}$, establishing $\gen{R_{\cF}}\subseteq k\pow{\cF}$.
		Conversely, let $A\in \cF$. Again by \Autoref{lem:rayner1}, we have $k\pow{\{A\}}\subseteq k\pow{\cF_A}=\gen{r_A}$.
		Hence,
		$$k\pow{\cF}=\bigcup_{A\in \cF}k\pow{\{A\}}\subseteq \bigcup_{A\in \cF}\gen{r_A}=\gen{R_{\cF}}.$$
	\end{proof}

	\begin{remark}\label{rmk:rayner3}
			In \cite[Conditions~2.1]{kks}, several properties of $\cF$ are established such that the $k$-hull $k\pow{\cF}$ forms a distinguished substructure of $k\pow{G}$. If $\cF$ contains $\{0\}$ and is closed under subsets as well as finite unions, then the assumption of \Autoref{prop:rayner2} on $\cF$ holds. Hence, for any such $\cF$, its $k$-hull $k\pow{\cF}$ is determined by $R_{\cF}$. In particular, every Rayner field is determined by linear recurrence relations. \qed
	\end{remark}

	\begin{corollary}
		Suppose that $k\neq \bbF_2$ and that $\{0\}\in \cF$. Suppose further that $k\pow{\cF}$ is a subgroup of $(k\pow{G},+)$. Then $k\pow{\cF}=\gen{R_{\cF}}$.
	\end{corollary}
	
	\begin{proof}
		This follows immediately from \cite[Proposition~3.4]{kks} and 	\Autoref{rmk:rayner3}.
	\end{proof}

		\begin{example}\label{cor:smth}
		\Autoref{prop:rayner2} and \Autoref{rmk:rayner3} 
		gives us several examples of distinguished subfields and other algebraic substructures of $k\pow{G}$ that are determined by linear recurrence relations:
		\begin{enumerate}[label = (\roman*)]
			\item\label{cor:smth:1}
			Recall that $\cF_{\fin}$ consists of all finite subsets of $G$. Thus, $\gen{R_{\cF_{\fin}}}=k[G]$. This shows that the group ring $k[G]$ is determined by linear recurrence relations.
			\item 
			If $\cF$ consists of all well-ordered subsets of $G^{\geq 0}$, then we obtain $\gen{R_\cF}=k\pow{G^{\geq 0}}$. Hence, the valuation ring of $k\pow{G}$ is determined by linear recurrence relations.
			\item 
			Let $\kappa$ be an uncountable cardinal. If $\cF$ consists of all well-ordered subsets of $G$ whose cardinality is strictly less than $\kappa$, then $\gen{R_\cF}=k\pow{G}_\kappa$.
			Likewise, if $\cF$ consists of all well-ordered subsets of $G^{\geq 0}$ whose cardinality is strictly less than $\kappa$, then $\gen{R_\cF}=k\pow{G^{\geq 0}}_\kappa$. See \cite{kuhlmann} for further details on $\kappa$-bounded power series.
			\item Note that the valuation ideal $k\pow{G^{>0}}$ of $k\pow{G}$ is not determined by linear recurrence relations for $G\neq \{0\}$, as any non-empty set that is determined by linear recurrence relations contains $k$ by \Autoref{cor:description}.  \qed
		\end{enumerate}
	\end{example}

				We now give an example to show that not every Hahn field that is determined by linear recurrence relations is also a Rayner field. For instance, \Autoref{prop:finitesets} shows that the minimal Hahn field $k(G)$ is always determined by linear recurrence relations. However, $ k(G) $ is never a Rayner field in the case $\Char(k)=0$ and $G\neq \{0\}$, as the following example shows.
				
	\begin{example}\label{ex:k(G)-not-k-hull}
			Suppose that $\Char(k)=0$ and $G\neq\{0\}$. Let $g\in G^{>0}$ and recall from \Autoref{ex:k(G)-not-k-hull0} that $e^{t^g}=\sum_{n=0}^\infty \tfrac{t^{ng}}{n!}\notin k(G)$. 
			The power series $ (1-t^g)^{-1}= \sum_{n=0}^\infty t^{ng} $ (cf.\ Neumann's Lemma \cite[page~211]{neumann}) has support $\Span_\N(g)$ and belongs to $ k(G) $. 
			Since $e^{t^g}$ has the same support as $(1-t^g)^{-1}$ but does not belong to $k(G)$, this shows that $k(G)$ is not a $k$-hull and thus not a Rayner field. \qed
	\end{example}

	In the last part of this section, we relate $\gen{R_\cF}$ to $\gen{S(\cF)}$, which was introduced in \Autoref{def:fseq}.	
	Recall that by \Autoref{prop:finitesets} we have $\gen{S(\cF_{\fin})}=k(G)$, and thus by \Autoref{cor:smth}~\autoref{cor:smth:1} $$\gen{S(\cF_{\fin})}=\ff(k[G])=\ff(k\pow{\cF_{\fin}})=\ff(\gen{R_{\cF_{\fin}}}).$$
	The following result shows that $\gen{S(\cF)}=\ff(\gen{R_{\cF}})$ also holds for a wider class of families $\cF$.
	
	\begin{proposition}\thlabel{condition-F-sequence}
		Suppose that $\cF$ contains $\{0\}$ and is closed under subsets, finite unions as well as addition. Then $\gen{S(\cF)}=\ff(\gen{R_{\cF}}) = \ff(k\pow{\cF})$.
	\end{proposition}

	\begin{proof}
		By \Autoref{cor:field}, $\gen{S(\cF)}$ is a field; by \Autoref{prop:rayner2} $\gen{R_{\cF}} = k\pow{\cF}$.
		For the inclusion $\gen{S(\cF)}\supseteq \ff(\gen{R_{\cF}})$, it thus suffices to show that $k\pow{\cF}\subseteq \gen{S(\cF)}$. To this end, let $a\in k\pow{G}$ with $\supp(a)\in \cF$. Then $r_{\supp(a)}$ has domain $\supp(a)\cup\{0\}$ and thus this sequence lies in $S(\cF)$. Hence,
		$$\frac{a}{r_{\supp(a)}^*}=\frac{a}{1}\in \gen{S(\cF)}.$$
		For the inclusion $\gen{S(\cF)}\subseteq  \ff(\gen{R_{\cF}})$, let $r\in S(\cF)$ and let $a\in k\pow{G}$ with $\supp(a)\subseteq \dom(r)\in \cF$. Since $\cF$ is closed under subsets, also $\supp(a),\allowbreak \supp(r)\allowbreak \in \cF$. Hence,
		$$\frac{a}{r^*}\in\ff(\gen{R_{\cF}}), $$
		as required.
	\end{proof}

	Since any Rayner field is a field (cf.\ \cite[Theorem~3.1]{kks}) and thus the fraction field of itself, \Autoref{condition-F-sequence} immediately implies the following.
	
	\begin{corollary}
		Suppose that $k\pow{\cF}$ is a Rayner field. Then $\gen{S(\cF)}=\gen{R_{\cF}}=k\pow{\cF}$.
	\end{corollary}

\section{Power series fields in finitely many variables}\label{sec:multivariate}

	This section exploits our results from the previous sections, in particular \Autoref{sec:hahnfields}, for Hahn fields that are fraction fields of formal power series rings in $n$ variables. We make precise in the following how these fields can be regarded as Hahn fields in the sense of \Autoref{def:hahn}.

	Let $ n\in\N\setminus\{0\} $. Throughout this section, we regard $\Z^n$ as an ordered abelian group with respect to the lexicographic ordering $<_{\lex}$.
	Denote by $ k[\x] = k[x_1,\ldots,x_n] $ the ring of polynomials in  the variables $ x_1,\ldots,x_n $. Let $ k(\x)=\ff (k[\x]) $ and let $ k\rpow{\x} $ be the ring of formal power series in the variables $ x_1,\ldots,x_n $ over $ k $, i.e.\ $k\rpow{\x}=k\rpow{x_1}\rpow{x_2}\ldots\rpow{x_n}$.
	Using multi-index notation, an element $ s\in k\rpow{\x} $ will be denoted by $ s = \sum_{g\in\N^n}s_{g}\x^{g} $, where $g=(g_1,\ldots,g_n)$, $ s_{g}\in k $ and $ \x^{g} = x_1^{g_1}\cdots x_n^{g_n} $.

	Since $\N^n$ is a well-ordered subset of $\Z^n$, the map $$  \sum_{g\in\N^n}s_{g}\x^{g}  \mapsto  \sum_{g\in\N^n}s_{g}t^{g} $$ is a $k$-isomorphism from $ k\rpow{\x} $ to $k\pow{\N^n} = \{s\in k\pow{\Z^n}\mid \supp(s)\subseteq \N^n\} $. We identify $ k\rpow{\x} $ with $k\pow{\N^n}$ via this $k$-isomorphism and thus regard $k\rpow{\x}$ as a subring of {the valuation ring of} $k\pow{\Z^n}$. Note that then $k[\x]= \{s\in k\pow{\Z^n}\mid \supp(s)\subseteq \N^n, |\supp(s)|<\infty\}$.
	Moreover, we write $ k\pow{\x}:= \ff (k\rpow{\x}) $ and obtain a $k$-embedding of $ k\pow{\x} $ into $ k\pow{\Z^n} $ via the identification above. This also enables us to expand an element $\tfrac p q\in k\pow{\x}$ with $p,q\in k\rpow{\x}$ as a power series of the form
	\begin{align*}\frac p q=\sum_{g\in\Z^n}s_{g}\x^{g},
	\end{align*}
	where $\sum_{g\in\Z^n}s_{g}t^{g}=\tfrac pq\in k\pow{\Z^n}$.
	
	\begin{remark}\thlabel{rmk:k((x))props}
		 
		\begin{enumerate}[label = (\roman*)]
		\item By the identifications described above, we obtain
		$$k(\Z^n)=k(\ul{x})\subseteq k\pow{\x}\subseteq k\pow{\Z^n}\!.$$
		Hence, both $k(\ul{x})$ and $k\pow{\x}$ are Hahn fields with coefficient field $k$ and value group $\Z^n$.
				
		\item\label{rmk:k((x))props:2} Since $\N^n$ is well-ordered, its power set $\mathcal{P}(\N^n)$ only consists of well-ordered subsets of $\Z^n$. Note that $k\rpow{\x}=k\pow{\N^n}=k\pow{\mathcal{P}(\N^n)}$. Since $\mathcal{P}(\N^n)$ satisfies the conditions of \Autoref{condition-F-sequence}, we obtain $$\gen{S(\mathcal{P}(\N^n))} = \ff(k\pow{\mathcal{P}(\N^n)})=k\pow{\x}\!.$$
		Hence, $k\pow{\x}$ is determined by linear recurrence relations. Since $\N^n$ is an upper bound of $\mathcal{P}(\N^n)$ with respect to the partial ordering $\subseteq$, \Autoref{prop:upperbound} shows that
		$$k\pow{\x}=\gen{S(\{\N^n\})},$$
		i.e.\ $k\pow{\x}$ is the subset of $k\pow{\Z^n}$ determined by $\{\N^n\}$-sequences.
		
		\item A cardinality computation shows that $k\pow{\x}$ is not the minimal Hahn field $k(\x)$: since $k\pow{\mathcal{P}(\N^n)}=k\rpow{\x}$, we obtain  \begin{align*}|k(\x)|&=\max\{|k|,\aleph_0\}<\max\setbr{2^{|k|},2^{\aleph_0}}\\&=(2^{\aleph_0})^{|k|}=|k\rpow{\x}|=|k\pow{\x}\!|.\end{align*}
		
		\item\label{rmk:k((x))props:4} Note that for the univariate case $n=1$ we obtain $k\pow{x_1}=k\pow{\Z}$, i.e.\ $k\pow{x_1}$ is already the maximal Hahn field. For $n\geq 2$, the field $k\pow{\ul{x}}$ is not the maximal Hahn field $k\pow{\Z^n}$, as witnessed by the following power series $s\in k\pow{\Z^n}\setminus k\pow{\x}$:
		Let $$s=\sum_{i=0}^\infty x_1^{i^2}x_2^{-i}=\sum_{i=0}^\infty t^{(i^2,-i,0,\ldots,0)}\in k\pow{\Z^n}\!.$$
		Note that $\supp(s)=\{(i^2,-i,0,\ldots,0)\mid i\in \N\}$ is well-ordered, as $$(0,0)<_{\lex}(1,-1)<_{\lex}(4,-2)<_{\lex}\ldots.$$
		Now assume, for a contradiction, that $s\in k\pow{\x}$. By \autoref{rmk:k((x))props:2} above, there is a (non-trivial) $\{\N^n\}$-sequence $r\colon \N^n\to k$  that determines $s$.
		Hence, for any $h=(h_1,\ldots,h_n)\in \Z^n\setminus \N^n$ the linear recurrence relation
		$$\sum_{g\in \N^n}r(g)s_{h-g}=0$$
		holds.
		By rewriting this sum, replacing $r$ by its associated power series $r^*\colon \Z^n\to k$ and removing $0$ terms, we obtain
		\begin{align}
			0&=\sum_{g\in \Z^n}r^*(g)s_{h-g}=\sum_{g_1,\ldots,g_n\in\notag \Z}r^*(g_1,\ldots,g_n)s_{(h_1-g_1,\ldots,h_n-g_n)}\notag\\
			&=\sum_{g_1,\ldots,g_n\in \Z}r^*(h_1-g_1,\ldots,h_n-g_n)s_{(g_1,\ldots,g_n)}\notag\\
			&=\sum_{i\in \N}r^*(h_1-i^2,h_2+i,h_3,\ldots,h_n)s_{(i^2,-i,0,\ldots,0)}\notag\\
			&=\sum_{i=0}^\infty r^*(h_1-i^2,h_2+i,h_3,\ldots,h_n).\label{eq:nozeroterms}
		\end{align}
		We show that $r^*=0$, leading to the required contradiction by the choice of $r$ as a non-trivial linear recurrence sequence. First note that for any $h_1,h_2\in \Z$ and any $(h_3,\ldots,h_n)\in \Z^{n-2}\setminus \N^{n-2}$ we have
		$r^*(h_1,\ldots,h_n)=0$ as $\supp(r^*)\subseteq \N^n$. We fix $h_3,\ldots,h_n\in \N^{n-2}$ and set $u(j,m)=r^*(j,m,h_3,\ldots,h_n)$ for any $j,m\in \Z$. It now suffices to show that for any $j,m\in \N$ we have $u(j,m)=0$. 	
		By (\autoref{eq:nozeroterms}), we have \begin{align}\sum_{i=0 }^\infty u(h_1-i^2,h_2+i)=0\label{eq:nozeroterms2}\end{align} for any $(h_1,h_2)\in \Z^2\setminus \N^2$.
		For $j,m\in \N$ we set $h_1=(m+j+1)^2+j\in \N$ and $h_2=-j-1\in \Z\setminus \N$.
		Since $\supp(r^*)\subseteq \N^n$ and $(m+j+1)^2+j-((m+1)+j+1)^2<0$, we obtain from (\autoref{eq:nozeroterms2}) that
		\begin{align}0&=\sum_{i=0 }^\infty u((m+j+1)^2+j-i^2,-j-1+i)\notag\\
			&= \sum_{i=0 }^\infty u((m+j+1)^2+j-(i+j+1)^2,i)\notag\\
			&= \sum_{i=0 }^m u((m+j+1)^2+j-(i+j+1)^2,i).\label{eq:nozeroterms3}
		\end{align}
		Using (\autoref{eq:nozeroterms3}), we show by induction on $m\in \N$ that for any $j,m\in \N$ we have $u(j,m)=0$.
		For $m=0$ we immediately obtain
		$$0=u((j+1)^2+j-(j+1)^2,0)=u(j,0).$$
		Now suppose that our induction hypothesis has been proved up to some $m\in \N$ and let $j\in \N$. Then
		\begin{align*}
			0&=\sum_{i=0 }^{m+1} u((m+j+2)^2+j-(i+j+1)^2,i)\\
			&= u((m+j+2)^2+j-((m+1)+j+1)^2,m+1)\\
			&=u(j,m+1),
		\end{align*}
		as required to complete the induction.		
		\end{enumerate}
	\end{remark}

	While \Autoref{rmk:k((x))props}~\autoref{rmk:k((x))props:4} above shows that $k\pow{\x}$ is not the maximal Hahn field for any $n\geq 2$, it leaves open whether $k\pow{\x}$ is a Rayner field. We pose this as an open question.\footnote{ Following a first preprint of our work, a negative answer to this question was presented in \cite[Example 2.11 and Remark 2.12]{hickel3}.}
	\begin{question}\thlabel{qu:rayner}
		Let $n\geq 2$. Is $k\pow{x_1,\ldots,x_n}$ a Rayner field?
	\end{question}
	In order to give a negative answer to \Autoref{qu:rayner}, it would suffice to find two power series $s,s'\in k\pow{\Z^n}$ with $\supp(s)=\supp(s')$, $s\in k\pow{\x}$ and $s'\notin k\pow{\x}$, as this would show that $k\pow{\x}$ is not a $k$-hull.
	
	We now address the query of identifying elements of $k(\x)$ within $k\pow{\x}$. A direct specialisation of  \Autoref{lem:coefficientskG} to $k\pow{\x}$ gives the following.
	\begin{corollary}\thlabel{theo-multivar} 
		Let $ s = \sum_{g\in \Z^n}s_g\x^g\in k\pow{\x} $. 
		Then $s\in k(\x)$ if and only if the following holds:
		there exist $\ell\in \N$, finite sets $A=\{h_0,\ldots,h_\ell\}\subseteq \Z^n$ with $h_0<_{\lex}\ldots<_{\lex}h_\ell$ and $\{r_0,\ldots,r_\ell\}\subseteq k$, not all $r_i$ equal to $0$, such that for any $g\in \Z^n\setminus A$ the following linear recurrence relation holds:
		\begin{align}
		\sum_{j=0}^\ell r_js_{g-h_j}=0\label{eq:lin-rec-multivar}
		\end{align}
	\end{corollary}

		For two non-zero polynomials $p,q\in k[\x]$, the sets $A=\{h_0,\ldots,h_\ell\}\subseteq \Z^n$ (we can even choose $A\subseteq \N^n$) and $\{r_0,\ldots,r_\ell\}\subseteq k$ such that the linear recurrence relation (\autoref{eq:lin-rec-multivar}) holds for $s=\tfrac pq$ are given by $A=\supp(p)\cup \supp(q)$ and $r_i=q_{h_i}$. The following example illustrates how this can be applied in practice in the bivariate case.
		
	\begin{example}
		Let $m\in \N$ and let $p,q\in k[x_1,x_2]\setminus \{0\}$ be given by $$p=\sum_{i=0}^m\sum_{j=0}^m p_{ij}x_1^ix_2^j \text{ and }q=\sum_{i=0}^m\sum_{j=0}^m q_{ij}x_1^ix_2^j.$$
		(Note that we do not require particular coefficients of $p$ or $q$ to be non-zero.)
		Consider $s=\tfrac pq\in k\pow{x_1,x_2}$. For any $i,j\in\{0,\ldots,m\}$, set $h_{ij}=(i,j)$ and $A = \{h_{ij}\mid i,j\in\{0,\ldots,m\}\}= \{0,\ldots,m\}^2\subseteq \N^2$. 
		Note that $$h_{00}<_{\lex}h_{01}<_{\lex}\ldots<_{\lex}h_{0m}<_{\lex}h_{10}<_{\lex}\ldots<_{\lex} h_{1m}<_{\lex}\ldots <_{\lex}h_{mm}.$$
		The linear recurrence sequence $r\colon \{0,\ldots,m\}^2\to k, (i,j)\mapsto q_{ij}$ has the property $r^*=q$. Moreover, $\supp(p)\subseteq \{0,\ldots,m\}^2=\dom(r)$. Hence, $s=\tfrac{p}{r^*}$ and by our Main Lemma (\Autoref{prop:description}) we obtain that for any $(g_1,g_2)\in \Z^2\setminus A$ the following linear recurrence relation holds:
		\begin{align*}
		\sum_{i=0}^m\sum_{j=0}^m q_{ij}s_{g_1-i,g_2-j}=0.
		\end{align*}\qed
	\end{example}
		
		As a converse to the above, given a power series $s\in k\pow{\x}$ satisfying (\autoref{eq:lin-rec-multivar}), we present in the following a procedure to construct $p,q\in k[\x]$ with $s=\tfrac pq$ by providing explicit formulas.
		
	\begin{construction}\thlabel{constr:pq}
		Let $ s = \sum_{g\in \Z^n}s_g\x^g\in k\pow{\x} $ and suppose that there exist $\ell\in \N$, $A=\{h_0,\ldots,h_\ell\}\subseteq \Z^n$ with $h_0<_{\lex}\ldots<_{\lex}h_\ell$ and $\{r_0,\ldots,r_\ell\}\subseteq k$, not all $r_i$ equal to $0$, such that for any $g\in \Z^n\setminus A$
		$$\sum_{j=0}^\ell r_js_{g-h_j}=0.$$ We compute the coefficients of $p,q\in k[\x]$ such that $s=\tfrac pq$.
		\begin{enumerate}[label = (\roman*)]
			\item\label{constr:pq:1} \emph{Computing $q$:} By \Autoref{prop:description}, we can express $s$ as
			$s=\tfrac{a}{b}$ for some $a,b\in k\pow{\x}$ with
			$b= \sum_{i=0}^\ell r_i\x^{h_i}$ and $\supp(a)\subseteq A$. However, despite their finite supports neither $a$ nor $b$ are necessarily polynomials in $k[\ul{x}]$, since the exponents $h_i\in \Z^n$ may contain negative entries. For any $i\in \{0,\ldots,\ell\}$, we write $h_i=( h_{i1},\ldots,h_{in})$. Let $\alpha=(\alpha_1,\ldots,\alpha_n)$ with $\alpha_j = -\min\{h_{ij} \mid i\in\{0,\ldots,\ell\} \}$ for any $j\in \{1,\ldots,n\}$. 
			Set \begin{align}q= \x^\alpha b= \sum_{i=0}^\ell r_i\x^{h_i+\alpha}\in k[\x],\label{eq:constr:pq:1}\end{align} as $h_i+\alpha\in \N^n$ for any $i\in \{0,\ldots,\ell\}$. 
			Then also $\x^\alpha a\in k[\x]$ and we obtain
			$$s=\frac {\x^\alpha a}{\x^\alpha b}\in k(\x).$$
			\item \emph{Computing $p$:} By step~\autoref{constr:pq:1}, there is $p=qs\in k[\x]$, where $q$ is given by the formula (\autoref{eq:constr:pq:1}). Moreover, the support of $p$ is given by $\supp(\x^\alpha a)$, which is contained in
			$\{h_i+\alpha\mid i\in \{0,\ldots,\ell\}\}$. Hence, we can express $p$ as
			$$p=\sum_{i=0}^\ell p_i\x^{h_i+\alpha}.$$
			To give an explicit formula for each $p_i$, consider the equation $p=qs$ amounting to
			$$\sum_{i=0}^\ell p_i\x^{h_i+\alpha}=\brackets{\sum_{j=0}^\ell r_j\x^{h_j+\alpha}}\cdot \brackets{\sum_{g\in \Z^n}s_g\x^g}\!.$$
			Thus, for any $i\in \{0,\ldots,\ell\}$ we obtain by comparing coefficients
			$$p_i=\sum_{j=0}^\ell \sum_{\substack{g\in \Z^n \\ h_i+\alpha=h_j+\alpha+g}}\!\!\!\!\!\!\!\!\!r_js_{g}=\sum_{j=0}^\ell r_js_{h_i-h_j}.$$
			This finally leads to the explicit formula for computing $p$:
			\begin{align}
				p=\sum_{i=0}^\ell \sum_{j=0}^\ell r_js_{h_i-h_j}\x^{h_i+\alpha}.\label{eq:constr:pq:2}
			\end{align}
		\end{enumerate} \qed
	\end{construction}

	\begin{remark} 
			Note that in \Autoref{constr:pq} the computation of $q$ and $p$ requires only finitely many coefficients of $s$: namely, the formula (\autoref{eq:constr:pq:1}) for $q$ is independent of the coefficients $s_g$ of $s$, and the formula (\autoref{eq:constr:pq:2}) for $p$ only requires the at most $(\ell+1)^2$ many coefficients $s_{h_i-h_j}$ for $(i,j)\in \{0,\ldots,\ell\}^2$. Hence, only finitely many ``initial values'' in the form of coefficients of $s$ together with the finite linear recurrence relation are needed to determine $p$ and $q$.\qed
	\end{remark}

	In the final part of this section, we present conditions on the linear recurrence relation (\autoref{eq:lin-rec-multivar}), and thus implicitly also on the corresponding polynomials $p\in k[\x]$ and $q\in k[\x]$ computable by \Autoref{constr:pq}, in order that $s=\tfrac pq$ already lies in $k\rpow{\x}$.
	
	\begin{corollary}\thlabel{prop:k[[x]]} 
		Let $ s = \sum_{g\in \Z^n}s_g\x^g\in k\pow{\x} $. 
		Set $h_0=0$ and let $r_0\in k\setminus\{0\}$.
		Suppose that
		there exist $\ell\in \N$, $A=\{h_0,\ldots,h_\ell\}\subseteq \N^n$ with $h_0<_{\lex}\ldots<_{\lex}h_\ell$ and $\{r_0,\ldots,r_\ell\}\subseteq k$ such that for any $g\in \Z^n\setminus A$
		\begin{align*}
		\sum_{j=0}^\ell r_js_{g-h_j}=0.
		\end{align*}
		Then $s\in k(\x)\cap k\rpow{\x}$.
	\end{corollary}

	\begin{proof} 
		\Autoref{theo-multivar} shows that $s\in k(\x)$. We have to verify that $\supp(s)\subseteq \N^n$. For some $p,q\in k[\x]$ of the form
		$$p=\sum_{i=0}^\ell p_ix^{h_i}\text{ and }q=r_0+\sum_{i=1}^\ell r_ix^{h_i},$$
		we have $s=\tfrac pq$. 
		Let $u_i=-r_0^{-1}r_i$ for $i\in \{0,\ldots,\ell\}$.
		We compute
		\begin{align*}
			s=\frac{\sum_{i=0}^\ell p_0x^{h_i}}{r_0\cdot\brackets{1-\sum_{i=1}^\ell u_ix^{h_i}}}
			=r_0^{-1}\sum_{i=0}^\ell p_0x^{h_i}\sum_{j=0}^\infty \brackets{\sum_{i=1}^\ell u_ix^{h_i}}^j\!\!.
		\end{align*}
		Hence,
		$$\supp(s)\subseteq A+\Span_{\N}(A)\subseteq \Span_{\N}(A)\subseteq \N^n,$$
		as $A\subseteq \N^n$. 
	\end{proof}

	\begin{example}
		 
		\begin{enumerate}[label = (\roman*)]
			\item 
			In the conditions on $A$ in \Autoref{prop:k[[x]]} it is not enough to require $A\subseteq \Z^n$ rather than $A\subseteq \N^n$, as the following example for $n=2$ shows.
			Consider $$s=\sum_{i=0}^\infty x_1^ix_2^{-i}=\frac{1}{1-x_1x_2^{-1}}=\frac{x_2}{x_2-x_1}\in k(x_1,x_2)\setminus k\rpow{x_1,x_2} $$
			with $\supp(s)=\{(i,-i)\mid i\in \N\}$. Setting $h_0=(0,0)$, $h_1=(1,-1)$, $r_0=1$ and $r_1=-1$, we obtain that $A=\{h_0,h_1\}\not\subseteq \N^2$ but for any $g\in \Z^2\setminus A$ the linear recurrence relation $s_g-s_{g-(1,-1)}=0$ holds.
			
			\item 			
			There are $s\in k(\x)\cap k\rpow{\x}$ for which there exist 
			$\ell\in \N$, $A=\{h_0,\ldots,h_\ell\}\allowbreak \subseteq \Z^n$ with $0=h_0<_{\lex}\ldots<_{\lex}h_\ell$ and $\{r_0,\ldots,r_\ell\}\subseteq k$ with $r_0\neq 0$ such that for any $g\in \Z^n\setminus A$
			\begin{align*}
			\sum_{j=0}^\ell r_js_{g-h_j}=0,
			\end{align*}
			but $A\not\subseteq \N^n$. For instance in the case $n=2$, choose $\ell=1$, $h_1=(1,-1)$ and $r_0=r_1=1$. Then $A=\{(0,0),(1,-1)\}\not\subseteq \N^2$ and $s=1$ can be expressed as
			$1 = \tfrac{1-x_1x_2^{-1}}{1-x_1x_2^{-1}}$ and thus satisfies the linear recurrence relation
			$\sum_{j=0}^\ell r_js_{g-h_j}=0$ for any $g\in \Z^2\setminus A$.			
			\qed
		\end{enumerate}
	\end{example}

		While \Autoref{prop:k[[x]]} gives \emph{sufficient} conditions on the coefficients of a power series $s\in k\pow{\x}$ in order that it lies in $k(\x)\cap k\rpow{\x}$, it remains open whether this condition is also \emph{necessary}. We pose this as a question.
	
	\begin{question}\thlabel{ass:conversek[[x]]}
		
		Let $ s = \sum_{g\in \Z^n}s_g\x^g\in k(\x)\cap k\rpow{\x} $. {Do there exist $\ell\in \N$, $A=\{h_0,\ldots,h_\ell\}\subseteq \Z^n$ with $h_0<_{\lex}\ldots<_{\lex}h_\ell$ and $\{r_0,\ldots,r_\ell\}\subseteq k$ (cf. \thref{theo-multivar}) that, additionally, satisfy $\{h_0,\ldots,h_\ell\}\subseteq \N^n$ and $r_0\neq0$, such that for any $g\in \Z^n\setminus A$ the linear recurrence relation
			\begin{align*}
				\sum_{j=0}^\ell r_js_{g-h_j}=0
			\end{align*}
			holds?}
	\end{question}

\section{Lifting properties}\label{sec:lifting}

	In this final section, we study Hahn fields that satisfy certain lifting properties as introduced in \cite{kuhlmann-serra-fields}. While there are several classes of lifting properties, we concentrate on the \emph{canonical first lifting property} (see \cite[Definition~3.5.5]{kuhlmann-serra-fields} and \Autoref{def:can1stlp} below) and the \emph{canonical second lifting property} (see \cite[Definition~3.3.5]{kuhlmann-serra-fields} and \Autoref{def:can2ndlp} below).
	
	\textbf{Throughout this section, let $K$ be a Hahn field. }
	Recall that $K$ is endowed with the canonical valuation $ v(s)=\min\supp(s) $.
	Let $\Aut(K)$ be the group of field automorphisms of $K$ and let $\oAut(G)$ be the group of order preserving group automorphisms of $G$. We say that $ \sigma \in \Aut(K)$ is \textbf{valuation preserving} if for any $ a,b\in K $, whenever $ v(a)=v(b)$, then $v(\sigma(a))=v(\sigma(b)) $.
	(For a characterisation of order preserving automorphisms, in terms of the valuation invariants of $ K $, see \cite[Theorem~4.2]{kuhlmannmatusinskipoint}.) The set of valuation preserving automorphisms of $ K $ forms a subgroup of $ \Aut(K) $ that we denote by $ \vAut(K) $.

\subsection{Canonical first lifting property}\label{sec:lifting1}
	 
	Let $ \sigma\in\vAut(K) $. We canonically associate to $\sigma $ the automorphism $ \sigma_k\in\Aut(k) $ given by $ \sigma_k (x) = (\sigma(x))_0 $ (where $\sigma(x)$ is regarded as an element in $k\pow{G}$), for any $ x\in k $, as well as the order preserving automorphism $ \sigma_G\in\oAut(G)$ given by $ \sigma_G(v(a)) = v(\sigma(a)) $, for any $ a\in K^\times $. This yields a homomorphism
	\[
	\Phi_K\colon \vAut (K) \to \Aut (k) \times \oAut (G),\ \sigma \mapsto (\sigma_k,\sigma_G).
	\]
	Notice that for the maximal Hahn field $ k\pow{G} $ the map $ \Phi_{k\,\pow{G}} $ admits the canonical section
	$\Psi_c\colon \Aut (k) \times \oAut (G) \to \vAut (k\pow{G})$ mapping $(\rho,\tau)\in  \Aut (k) \times \oAut (G)$ to 
	$$
	\Psi_c(\rho,\tau)\colon  k\pow{G} \to k\pow{G},
	 \sum_{g\in G} s_g t^g \mapsto  \sum_{g\in G} \rho(s_g)t^{\tau(g)}.$$
	Following the notation of \cite[Definition~3.5.2]{kuhlmann-serra-fields}, we also denote $\Psi_c(\rho,\tau)$ simply by $\widetilde{\rho\tau}$ and call it the \textbf{canonical lift} of $(\rho,\tau)$.
	
	\begin{definition}\thlabel{def:can1stlp}
		We say that $ K $ has the \textbf{canonical first lifting property} if for any $ (\rho,\tau)\in \Aut (k) \times \oAut (G) $ we have $ \widetilde{\rho\tau}|_K \in \vAut(K) $, i.e.\ if $ \Psi_{c,K}\colon (\rho,\tau) \mapsto \Psi_c(\rho,\tau)|_K $ is a section of $ \Phi_K $.
	\end{definition}

	Note that $K$ has the canonical first lifting property if and only if for any $ (\rho,\tau)\in \Aut (k) \times \oAut (G) $ and any $s= \sum_{g\in G} s_g t^g\in K$, also $\widetilde{\rho\tau}(s)=\sum_{g\in G} \rho(s_g)t^{\tau(g)}$ lies in $K$. 
	Note also that for any $ s\in K $ we have $ \supp(\widetilde{\rho\tau}(s)) \allowbreak = \tau(\supp(s)) $.
	For a general and detailed study of the canonical first lifting property (and others) we refer the reader to \cite{kuhlmann-serra-fields, serra-thesis}.
	
	\Autoref{not:basicnotions} introduced a way to associate to a power series $s\in k\pow{G}$ a linear recurrence sequence $s^*\in \LRS(k,G)$. In the following, we establish a similar notion for maps on $k\pow{G}$, i.e.\ we associate to a given $\varphi\colon k\pow{G}\to k\pow{G}$ a map $\varphi^*\colon \LRS(k,G)\to\LRS(k,G)$.
	{  To this purpose, we start with $r\in \LRS(k,G)$. We obtain from it an auxiliary $ 0 $-free $\widetilde r\in \LRS(k,G)$ satisfying $\dom(\widetilde{r})=\supp(\widetilde{r})=\dom(r)$, which we use to define the domain of $ \phi^*(r) $. While the definition of $ \phi^* $ might look artificial as it shows a strong dependence on the choice of $ \widetilde r $, we will see in \Autoref{lem:conliftdescr} that if $ \phi $ is a canonical lift, then this dependence disappears.}
		
	\begin{definition}
		\begin{enumerate}[label = (\roman*)]
			\item 	Let $r\in \LRS(k,G)$. We set $\widetilde{r}$ to be the linear recurrence sequence with domain $\dom(r)$ and
			$$\widetilde{r}(g):=\begin{cases}
			r(g)& \text{ for } g\in\supp(r),\\
			1& \text{ for } g\in\dom(r)\setminus \supp(r).
			\end{cases}$$
			
			\item Let $\varphi\colon k\pow{G}\to k\pow{G}$. We define the \textbf{map $\varphi^*\colon \LRS(k,G)\mapsto \LRS(k,G)$ associated} to $\varphi$ as follows: 
				For any $r\in \LRS(k,G)$ we let $$\dom(\varphi^*(r)):= \supp(\varphi(\widetilde{r}^*))$$ 
				and set
				$\varphi^*(r)(g):=
				\varphi(r^*)(g)$ for any $g\in \dom(\varphi^*(r))$.
		\end{enumerate}
	\end{definition}

		Note that for any $r\in \LRS(k,G)$ we have $\dom(\widetilde{r})=\supp(\widetilde{r})=\dom(r)$ and $\widetilde{r}^{**}=\widetilde{r}$. In particular, $\widetilde{r}$ is $0$-free. Moreover, if $r$ is already $0$-free, then $\widetilde{r}=r$. 
		Our choice of the map $\varphi^*\colon \LRS(k,G)\to \LRS(k,G)$ associated to $\varphi\colon k\pow{G}\to k\pow{G}$ makes the following diagram commute:
		\begin{center} 
			\begin{tikzcd}
				k\pow{G} \arrow{r}{\phi} \arrow{d}{*} &k\pow{G} \arrow{d}{*}\\
				\LRS(k,G) \arrow{r}{\phi^*} &\LRS(k,G)
			\end{tikzcd}
		\end{center}
		Indeed, let $s\in k\pow{G}$. Then $\dom(\varphi^*(s^*))=\supp(\varphi(s^{**}))=\supp(\varphi(s))=\dom(\varphi(s)^*)$, as $s^*$ is $0$-free. Moreover, for any $g\in \supp(\varphi(s))$, we have
		$\varphi^*(s^*)(g)= \varphi(s^{**})(g) =\varphi(s)(g)=\varphi(s)^*(g)$, establishing that $\varphi^*(s^*)=\varphi(s)^*$, as required.\footnote{
		Note that the above diagram suggests obtaining an algebraic structure on $\LRS(k,G)$ from that on $k\pow{G}$ via the maps of the diagram, which is beyond the scope of the present paper.
		}

		The following lemma gives a more explicit description of the map associated to a canonical lift of a pair $(\rho,\tau)\in \vAut(k)\times \oAut(G)$.
		
	\begin{lemma}\thlabel{lem:conliftdescr} 
		Let $(\rho,\tau)\in \vAut(k)\times \oAut(G)$ and consider $\varphi=\widetilde{\rho\tau}$. Let $r=(g_i,r_i)_{i<\alpha}\in \LRS(k,G)$. Then
		$\varphi^*(r)=(\tau(g_i),\rho(r_i))_{i<\alpha}$.
	\end{lemma}

	\begin{proof} 
		First we verify that the domain of $\varphi^*(r)$ is given by $\{\tau(g_i)\mid i<\alpha\}=\tau(\dom(r))$:
		$$
		\dom(\varphi^*(r))= \supp(\varphi(\widetilde{r}^*)) = \tau(\supp(\widetilde{r}^*))=\tau(\dom(r)).$$
		Since $\varphi(r^*)=\sum_{i<\alpha}\rho(r_i)t^{\tau(g_i)}$, we obtain
		for any $i<\alpha$ that $\varphi^*(r)(\tau(g_i))=\varphi(r^*)(\tau(g_i))=\rho(r_i)$, as required.
	\end{proof}

	In order to illustrate \Autoref{lem:conliftdescr}, we apply it to particular valuation preserving automorphisms and linear recurrence sequences to obtain specific examples.
		
	\begin{example} \thlabel{ex:dualitylift}
		\begin{enumerate}[label = (\roman*)]
			\item 	Since $\id_{k\,\pow{G}}$ is the canonical lift of $(\id_k,\id_G)$, we have $\id_{k\,\pow{G}}^*(r)=r$ for any $r\in \LRS(k,G)$, i.e.\ $\id_{k\,\pow{G}}^*=\id_{\LRS(k,G)}$.
			
			\item\label{ex:dualitylift:2} Let $A\subseteq G$ be well-ordered
			 and recall from \Autoref{not:rayner} that $r_A$ denotes the linear recurrence sequence with domain $A\cup\{0\}$ given by $r_A(0)=1$ and $r_A(g)=0$ for any $g\in A\setminus\{0\}$.
			 Let $(\rho,\tau)\in \Aut(k)\times \oAut(G)$ and $\varphi=\widetilde{\rho\tau}$. 
			 Then $\varphi^*(r_A)$ is the linear recurrence sequence with domain $\tau(A)\cup\{\tau(0)\}=\tau(A)\cup\{0\}$ given by $\varphi^*(r_A)(0)=\rho(1)=1$ and $\varphi^*(r_A)(g)=\rho(0)=0$ for any $g\in \tau(A)\setminus\{0\}$. Hence, $\varphi^*(r_A)=r_{\tau(A)}$. In particular, $\varphi^*(r_A)$ only depends on the automorphism $\tau$ on the value group and is independent of the automorphism $\rho$ on the residue field.
		\end{enumerate}
		\qed
	\end{example}

		We now establish a sufficient condition on a Hahn field determined by linear recurrence relations to have the canonical first lifting property.
		
		\begin{lemma}\thlabel{prop:phi*lift}
			Let $\varphi\in \Aut(k\pow{G})$ and suppose that
			for any $a,b\in k\pow{G}$ we have
			$$\supp(a)\subseteq \supp(b) \text{ if and only if } \supp(\varphi(a))\subseteq\supp(\varphi(b)).$$
			Then for any $r\in \LRS(k,G)$, 
			$$\varphi(\gen{r})=\gen{\varphi^*(r)}.$$
		\end{lemma}
	
	\begin{proof}
		If $r$ is trivial, then $\varphi^*(r)$ is also trivial, as $\varphi^*(r)(g)=\varphi(r^*)(g)=\varphi(0)(g)=0$ for any $g\in \dom(\varphi^*(r))$. Hence, $\varphi(\gen{r})=\varphi(k\pow{G})=k\pow{G}=\gen{\varphi^*(r)}$, as required.
		Suppose that $r$ is non-trivial.
		 
		\textbf{Claim:} $(\varphi^*(r))^*=\varphi(r^*)$.
		
		\noindent By definition of $\varphi^*$ we obtain
		\begin{align*}
			\supp(\varphi^*(r))& = \dom(\varphi^*(r))\cap \supp(\varphi(r^*))\\
			& = \supp(\varphi(\widetilde{r}^*)) \cap \supp(\varphi(r^*)).
		\end{align*}
		Since $\supp(\widetilde{r}^*)\supseteq \supp(r^*)$, we obtain $\supp(\varphi(\widetilde{r}^*))\supseteq \supp(\varphi(r^*))$ and thus
		$$\supp((\varphi^*(r))^*) = \supp(\varphi^*(r)) = \supp(\varphi(r^*)).$$
		Now for any $g\in \supp(\varphi(r^*))$, we have
		$$(\varphi^*(r))^*(g)=\varphi^*(r)(g)=\varphi(r^*)(g),$$
		establishing the claim.
		
		Let $a\in k\pow{G}$ with $\supp(a)\subseteq \dom(r)$.
		In order to show $\varphi(\gen{r})\subseteq\gen{\varphi^*(r)}$, it suffices to verify that $\supp(\varphi(a))\subseteq \dom(\varphi^*(r))$, as
		$$\varphi\!\brackets{\tfrac{a}{r^*}}=\tfrac{\varphi(a)}{(\varphi^*(r))^*}$$
		by the claim. 
		Since $\supp(a)\subseteq \dom(r)=\supp(\widetilde{r}^*)$, we obtain $$\supp(\varphi(a))\subseteq \supp(\varphi(\widetilde{r}^*))=\dom(\varphi^*(r)),$$
		as required.
		
		The converse $\varphi(\gen{r})\supseteq\gen{\varphi^*(r)}$ follows by similar arguments: If $b\in k\pow{G}$ with $\supp(b)\subseteq \dom(\varphi^*(r))=\supp(\varphi(\widetilde{r}^*))$, then 
		$$\supp(\varphi^{-1}(b))\subseteq \supp(\widetilde{r}^*)=\dom(r).$$
		Hence, for $a=\varphi^{-1}(b)$ we obtain
		$$\tfrac{b}{(\varphi^*(r))^*}=\tfrac{\varphi(a)}{\varphi(r^*)}=\varphi\!\brackets{\tfrac a{r^*}}\in \varphi(\gen{r}).$$
	\end{proof}

	\begin{theorem}\thlabel{thm:1lpcriterion} 
		Let $R\subseteq \LRS(k,G)$ such that $K=\gen{R}$ is a Hahn field. Suppose that for any $(\rho,\tau)\in \Aut(k)\times \oAut(G)$ we have $\widetilde{\rho\tau}^*(R)\subseteq R$.
		Then $K$ has the canonical first lifting property. 
	\end{theorem}

	\begin{proof} 
		Let $(\rho,\tau)\in \Aut(k)\times \oAut(G)$ and set $\varphi=\widetilde{\rho\tau}$. 
		We have to verify that $\varphi(K)\subseteq K$.
		
		Let $a,b\in k\pow{G}$. Since $\tau$ is bijective, we obtain 
		\begin{align*}
		\supp(a)\subseteq \supp(b)  \Leftrightarrow \underbrace{\tau(\supp(a))}_{=\supp(\varphi(a))}\subseteq \underbrace{\tau(\supp(b))}_{=\supp(\varphi(b))}.
		\end{align*}
		 We may thus apply \Autoref{prop:phi*lift} to $\varphi$, i.e.\ for any $r\in \LRS(k,G)$ we have	
		$\varphi(\gen{r})=\gen{\varphi^*(r)}$. We obtain
		\begin{align*}
		\varphi(K)=\varphi(\gen{R})=\bigcup_{r\in R}\varphi(\gen{r})
		=\bigcup_{r\in R}\gen{\varphi^*(r)}=\gen{\varphi^*(R)}\subseteq \gen{R}=K,
		\end{align*}
		as required.
	\end{proof}

	In the following, we apply \Autoref{thm:1lpcriterion} to Hahn fields determined by $\cF$-sequences (see \Autoref{sec:hahnfields}) as well as to $k$-hulls (see \Autoref{sec:rayner}).
	
	\begin{definition} 
		We say that $\cF$ is \textbf{stable under} $\oAut(G)$ if for any $\tau\in \oAut(G)$ and any $A\in \cF$ we have $\tau(A)\in \cF$.
	\end{definition}

	\begin{example}\thlabel{ex:lut}
		  Let $n\in \N\setminus\{0\}$
		and consider the group $ \Z^n $ ordered lexicographically by $<_{\lex}$. 
		Let $ \mathrm{LUT}_n(\Z) $ be the multiplicative group of lower unitriangular matrices with entries in $ \Z $, i.e. $ \mathrm{LUT}_n(\Z) $ consists of all 		$A=(a_{ij})_{i,j=1,\ldots,n}\in \Z^{n\times n}$ with $a_{ii}=1$ for any $i\in \{1,\ldots,n\}$ and $a_{ij}=0$ for any $i,j\in \{1,\ldots,n\}$ with $i<j$.
		By \cite[Corollary to Lemma~1]{conrad} we have $(\mathrm{LUT}_n(\Z), \cdot )\simeq (\oAut(\Z^n), \circ)  $ via the isomorphism $A \mapsto (g\mapsto Ag)$. Thus, for any $ \tau\in\oAut(\Z^n) $ there is a unique matrix $ A\in \mathrm{LUT}_n(\Z) $ such that for any $ \underline{m}=(m_1,\ldots,m_n)\in \Z^n $ we have\footnote{For details on this, see \cite[Section~2.3.2]{serra-thesis}. Notice that, there, the author works with upper unitriangular matrices, multiplying row vectors on the left.}
		\[
		\tau(\underline{m}) = (A\underline{m}^T)^T.
		\]
		A family $ \cF $ of subsets of $ \Z^n $ is therefore stable under $ \oAut(\Z^n) $ if and only if it is stable under the action of $ \mathrm{LUT}_n(\Z) $, i.e.\ for any $\Delta\in \cF$ and any $A\in\mathrm{LUT}_n(\Z)$ also $A\Delta = \{(A\underline{m}^T)^T\mid \underline{m}\in \Delta\}\in \cF$. 
		\qed
	\end{example}
	
	Recall from \Autoref{def:fseq} that an $\cF$-sequence is a non-trivial linear recurrence sequence whose domain lies in $\cF$, and the set of all $\cF$-sequences is denoted by $S(\cF)$.

	\begin{corollary}\thlabel{thm:canonicallifting} 
		Suppose that $\cF$ is stable under $\oAut(G)$. Then the following hold.
		\begin{enumerate}[label = (\roman*)]
			\item\label{thm:canonicallifting:1} If $\gen{S(\cF)}$ is a Hahn field, then it has the canonical first lifting property.
			\item\label{thm:canonicallifting:2} If $\cF_{\fin}\subseteq \cF$ and $\cF$ is closed under addition, then $\gen{S(\cF)}$ is a Hahn field with the canonical first lifting property.
		\end{enumerate}
	\end{corollary}

	\begin{proof} 
		By \Autoref{cor:hahnfield}, if $\cF_{\fin}\subseteq \cF$ and $\cF$ is closed under addition, then $K=\gen{S(\cF)}$ is already a Hahn field. We thus suppose that $\gen{S(\cF)}$ is a Hahn field and verify that $K$ has the canonical first lifting property.
		
		Let $(\rho,\tau)\in \Aut(k)\times \oAut(G)$ and set $\varphi=\widetilde{\rho\tau}$. Moreover, let $A\in \cF$ and let $r=(g_i,r_i)_{i<\alpha}\in \LRS(k,G)$ be non-trivial with $\dom(r)=A$.
		By \Autoref{thm:1lpcriterion} it now suffices to verify that $\varphi^*(r)\in S(\cF)$.
		\Autoref{lem:conliftdescr} shows that  $\varphi^*(r)=(\tau(g_i),\rho(r_i))_{i<\alpha}$. Thus, $\varphi^*(r)$ is non-trivial and its domain is given by $\tau(A)\in \cF$, i.e.\ $\varphi^*(r)$ is an $\cF$-sequence, as required.
	\end{proof}

	The following example shows that there are Hahn fields determined by $\cF$-sequences without the canonical first lifting property.
	
	\begin{example}\thlabel{ex:810}
		Let $n\in \N\setminus \{0,1\}$ and recall from \Autoref{rmk:k((x))props}~\autoref{rmk:k((x))props:2} that $k\pow{\x}=k\pow{x_1,\ldots,x_n}$ is determined by $S(\{\N^n\})$. 
		Let $A=I_n-E_{21}$,
		where $I_n\in \Z^{n\times n}$ denotes the identity matrix and $E_{21}=(\varepsilon_{ij})_{i,j=1,\ldots,n}\in \Z^{n\times n}$ is given by $\varepsilon_{21}=1$ and $\varepsilon_{ij}=0$ for any $(i,j)\neq (2,1)$.
		Then $(A(m_1,\ldots,m_n)^T)^T= (m_1,m_2-m_1,m_3,\ldots,m_n)$.
		Hence by \Autoref{ex:lut}, the map
		$$\tau\colon \Z^n \to \Z^n, (m_1,\ldots,m_n)\mapsto (m_1,m_2-m_1,m_3,\ldots,m_n)$$
		lies in $\oAut(\Z^n)$. Now $\tau(\N^n)\neq \N^n$ (e.g.\ $\tau(1,0,0,\ldots,0)= (1,-1,0,\ldots,0)\notin \N^n$), so $\{\N^n\}$ is not stable under $\oAut(\Z^n)$.
		
		Consider the power series $$s=\sum_{i=0}^\infty x_1^{i^2}x_2^{i^2-i}\in k\pow{\x},$$
		whose support is given by $\supp(s)=\{(i^2,i^2-i,0,\ldots,0)\mid i\in \N \}\subseteq \N^n$. Then $$\widetilde{\id_k\tau}(s)=\sum_{i=0}^\infty x_1^{i^2}x_2^{-i}\notin k\pow{\x}$$
		by \Autoref{rmk:k((x))props}~\autoref{rmk:k((x))props:4}. Hence, $k\pow{\ul{x}}$ does not have the canonical first lifting property.
	\end{example}

	\begin{corollary}\thlabel{cor:canonicallifting2} 
		Let $k\neq \mathbb{F}_2$ and suppose that $k\pow{\cF}$ is a field. Then it is a Hahn field with the canonical first lifting property if and only if $\{\{g\}\mid g\in G\}\subseteq \cF$ and $\cF$ is stable under $\oAut(G)$.
	\end{corollary}
	
	\begin{proof} 
		Note that $k\pow{\cF}$ is a Hahn field if and only if $\{t^g\mid g\in G\}\subseteq k\pow{\cF}$, and the latter holds if and only if $\{\{g\}\mid g\in G\}\subseteq \cF$.
		Since $k\pow{\cF}$ is a field and $k\neq \mathbb{F}_2$, \cite[Proposition~3.4]{kks} implies that $\cF$ is closed under subsets. Moreover, for any $A\in \cF$ with $0\notin A$, we have
		$s=\sum_{g\in A}t^g+1\in k\pow{\cF}$ and thus $A\cup\{0\}=\supp(s)\in \cF$. Hence, $\cF$ is also closed under unions with $\{0\}$. We can thus apply \Autoref{prop:rayner2} to obtain
		$\gen{R_{\cF}}=k\pow{\cF}$.
		
		Suppose that $\{\{g\}\mid g\in G\}\subseteq \cF$ and $\cF$ is stable under $\oAut(G)$. Let $(\rho,\tau)\in \Aut(k)\times \oAut(G)$ and set $\varphi=\widetilde{\rho\tau}$. Moreover, let $A\in \cF$.
		As shown in 
		\Autoref{ex:dualitylift}~\autoref{ex:dualitylift:2}, we have $\varphi^*(r_A)=r_{\tau(A)}\in R_{\cF}$, as $\tau(A)\in \cF$. \Autoref{thm:1lpcriterion} now yields that $k\pow{\cF}$ has the canonical first lifting property.
		
		Conversely, suppose that $k\pow{\cF}$ is a Hahn field with the canonical first lifting property. Let $\tau\in \oAut(G)$ and let $A\in \cF$. Since $a=\sum_{g\in A}t^g\in k\pow{\cF}$, also $a'=\widetilde{\id_k \tau}(a)=\sum_{g\in A}t^{\tau(g)}\in k\pow{\cF}$. This shows that $\tau(A)=\supp(a')\in \cF$, as required.
	\end{proof}

	Note that \Autoref{cor:canonicallifting2} applies, in particular, to Rayner fields. This was explicitly proved in \cite[Proposition~3.6.6]{kuhlmann-serra-fields} (also for the case $k=\mathbb{F}_2$).
	The results above can now be used to construct Hahn fields with the canonical first lifting property that are determined by linear recurrence relations, as follows.
		
	\begin{construction} \thlabel{constr:1clpfields}
		Let $\cF$ be an arbitrary family of well-ordered subsets of $G$ and set $\cF_{-1}=\cF\cup \cF_{\fin}$.
		\begin{enumerate}[label = (\roman*)]
			\item For any even $i\in \N$, set
			$$\cF_i=\{A+B\mid A,B\in \cF_{i-1}\}\text{ and }\cF_{i+1}=\{\tau(A)\mid A\in \cF_{i},\tau\in \oAut(G)\}.$$
			Then $\cF'=\bigcup_{i=0}^\infty \cF_i$ is closed under addition and stable under \allowbreak  $\oAut(G)$. 
			Hence, \Autoref{thm:canonicallifting}~\autoref{thm:canonicallifting:2} shows that $\gen{S(\cF')}$ is a Hahn field with the canonical first lifting property.
			
			\item For any $i\in \N$ that is an integer multiple of $5$, set
			\begin{align*}
				\cF_i&=\{A\mid A\subseteq B\in \cF_{i-1}\},\\ \cF_{i+1}&=\{A\cup B\mid A,B\in \cF_{i}\},\\
				\cF_{i+2}&=\{A+g\mid A\in  \cF_{i+1},g\in G\},\\ \cF_{i+3}&={\cF_{i+2}\cup}\{\Span_{\N}(A)\mid  A\in \cF_{i+2}, A\subseteq G^{\geq 0}\},\\
				\cF_{i+4}&=\{\tau(A)\mid A\in \cF_{i+3},\tau\in \oAut(G)\}.
			\end{align*}
			Moreover, let $\cF'=\bigcup_{i=0}^\infty \cF_i$.
			Then $k\pow{\cF'}$ is the smallest Hahn field that is a Rayner field extending $k\pow{\cF}$ with its family of supports $\cF'$ being stable under $\oAut(G)$. Recall that by \Autoref{prop:rayner2}, $k\pow{\cF'}=\gen{R_{\cF'}}$. Finally, 
			\cite[Proposition~3.6.6]{kuhlmann-serra-fields} shows that $k\pow{\cF'}$ has the canonical first lifting property.
		\end{enumerate}\qed
	\end{construction}
	
	While \Autoref{constr:1clpfields} gives rise to classes of Hahn fields with the canonical first lifting property that are determined by linear recurrence relations, not every Hahn field determined by linear recurrence relations has the canonical first lifting property. In \Autoref{ex:raynernolift} below, we present an explicit Rayner field -- which is determined by linear recurrence relations by \Autoref{rmk:rayner3} -- that does not have the canonical first lifting property.
	
	\begin{lemma}\label{lemma:minimal-F} 
		Let $ A_1,\ldots,A_m $ be finitely many subsets of $ G $. {Moreover, let $ \cF $ consist of all well-ordered $B\subseteq G$ for which there exist $n\in \N$ and $g_1,\ldots,g_n\in G$ with $B\subseteq \Span_{\Z}(g_1,\ldots,g_n,A_1,\ldots,\allowbreak A_m)$.} Then $k\pow{\cF}$ is a Hahn field that is a Rayner field.
	\end{lemma}

	\begin{proof} 
		We verify the conditions on $\cF$ establishing that $k\pow{\cF}$ is a Rayner field.
		Clearly, $\cF$ is non-empty and closed under subsets. Moreover, \linebreak $\Span_\Z\!\brackets{\bigcup_{A\in \cF}A}=G$, as $\{\{g\}\mid g\in G\}\subseteq \cF$, which also implies that $k\pow{G}$ is a Hahn field.
		
		Let $B,C,D\in \cF$ with $D\subseteq G^{\geq 0}$ and let $h\in G$. Then there are $n\in \N$ and $g_1,\ldots,g_n\in G$ such that
		$$B,C,D\subseteq \Span_{\Z}(g_1,\ldots,g_n,h,A_1,\ldots,A_m)=:J.$$
		Now also $B\cup C,B+h\subseteq J$, showing that $\cF$ is closed under finite unions and translations. Finally, also $\Span_{\N}(D)\subseteq \Span_{\Z}(D) \subseteq J$, whence $\Span_{\N}(D)\in \cF$.
	\end{proof}

	\begin{example}\thlabel{ex:raynernolift}
		We denote by $G=\coprod_\Z \Q$ the symmetric Hahn sum over $(\Q,\Z)$, that is, $G$ is the additive group consisting of all sums of the form $\sum_{i=-m}^m a_i\one_i$ for some $m\in \N$ and $a_i\in \Q$, where $\one_i\colon \Z\to\Q$ is the indicator function on $\{i\}$ mapping $i$ to $1$ and everything else to $0$. The ordering on $G$ making it an ordered abelian group is given by $\sum_{i=-m}^m a_i\one_i<\sum_{i=-m}^m b_i\one_i$ if and only if $a_{j}<b_j$ for $j=\min\{i\in \{-m,\ldots,m\}\mid a_i\neq 0\text{ or } b_i\neq 0\}$.
		
		Let $A_1=\{-\tfrac1{p_i}\one_1\mid i\in \N\}$, where $p_i\in \N$ denotes the $(i+1)$-th prime number. Let $\cF$ consist of all well-ordered subsets of subgroups of $G$ of the form $  \Span_{\Z}( g_1,\ldots,g_n,A_1 ) $, where $n\in \N$ and $g_1,\ldots,g_n\in G$. By   \Autoref{lemma:minimal-F}, $\R\pow{\cF}$ is a Hahn field that is a Rayner field. Since $A_1$ is well-ordered, we have $A_1\in \cF$. Consider the automorphism $\tau$ on $G$ given by $$\tau\colon \sum_{i=-m}^m a_i\one_i\mapsto \sum_{i=-m}^m a_i\one_{i-1}.$$ Then $\tau(A_1)=\{-\tfrac1{p_i}\one_0\mid i\in \N\}$. We show that $\tau(A_1)$ cannot be a subset of a subgroup of $G$ of the form  $  \Span_{\Z}( g_1,\ldots,g_n,A_1 ) $, whence $\tau(A_1)\notin \cF$. \Autoref{cor:canonicallifting2} will then imply that $\R\pow{\cF}$ does not have the canonical first lifting property.
		
		To this end, let $n\in \N$ and $g_1,\ldots,g_n \in G$ be arbitrary and let $h\in  \Span_{\Z}( g_1,\ldots,g_n,A_1 ) $. Then there exist $\ell\in \N$, $a_1,\ldots,a_\ell \in A_1$ and $r_1,\ldots,r_n,\allowbreak s_1,\ldots,s_\ell\in \Z$ such that
		$$h=r_1g_1+\ldots+r_ng_n+s_1a_1+\ldots+s_\ell a_\ell.$$
		Evaluating at $0$, we obtain
		$$h(0)=r_1g_1(0)+\ldots+r_ng_n(0),$$
		as $a_i(0)=0$ for any $i\in \{1,\ldots,\ell\}$.
		Let $q_1,\ldots,q_m\in \N$ be all distinct primes dividing at least one of the denominators of $g_1(0),\ldots,g_n(0)\in \Q$ expressed in their reduced fraction form, and let $p\in \N$ be a prime number such that $p>q_i$ for all $i\in \{1,\ldots m\}$. 
		Then $g=-\tfrac{1}{p}\one_0\in \tau(A_1)$. Hence, $g(0)=-\tfrac 1 p$. However, the denominator of $r_1g_1(0)+\ldots+r_ng_n(0)$ in reduced fraction form is not divisible by $p$. Hence, $g(0)$ cannot be expressed as a $\Z$-linear combination of $g_1(0),\ldots,g_n(0)$. This shows that $g\notin   \Span_{\Z}( g_1,\ldots,g_n,A_1 ) $, as required. \qed
	\end{example}

\subsection{Canonical second lifting property}\label{sec:lifting2}

	
	Recall that the map $\Phi_K$ canonically associates to every $\sigma\in \vAut(K)$ a pair $(\sigma_k,\sigma_G)\in \Aut(k)\times\oAut(G)$. The kernel of $\Phi_K$ is a normal subgroup of $\vAut(K)$ consisting of all valuation preserving automorphisms on $K$ whose associated automorphisms on $k$ and $G$ are the respective identity maps. We denote by $\Int\Aut(K)$ the kernel of $\Phi_K$ and call an automorphism $\sigma\in \Int\Aut(K)$ an \textbf{internal automorphism} (see \cite[Section~3]{kuhlmann-serra-fields} for further reference). The canonical second lifting property regards a particular class of internal automorphisms that are induced by group homomorphisms from $G$ to $k^\times$.
	Given such a homomorphism $x\in \Hom(G,k^\times)$, we denote its image $x(g)$ by $x^g$ for any $g\in G$. 
	Consider the map
	$$P\colon \Hom(G,k^\times)\to \Int \Aut(k\pow{G}), x\mapsto \rho_x,$$
	where $\rho_x$ is given by
	$$\rho_x\colon k\pow{G}\to k\pow{G}, \sum_{g\in G}s_gt^g \mapsto \sum_{g\in G}s_gx^gt^g.$$

	\begin{definition}\thlabel{def:can2ndlp} 
		We say that $ K $ has the \textbf{canonical second lifting property} if for any $ x\in \Hom(G,k^\times) $ we have $ \rho_x|_K \in \Int\Aut(K) $, i.e.\ if $ P_{K}\colon x \mapsto P(x)|_K $ is well-defined group homomorphism from $\Hom(G,k^\times)$ to $\Int\Aut(K)$.
	\end{definition}
	
	The map $ P\colon \Hom(G,k^\times) \rightarrow\Int\Aut( k\pow{G} )$ is always a well-defined group homomorphism, which implies that the maximal Hahn field $ k\pow{G} $ always has the canonical second lifting property.
	For any $ x\in \Hom(G,k^\times) $ and any $ s\in k\pow{G} $ we have $ \supp(s) = \supp(\rho_x(s)) $.
	Note that $K$ has the canonical second lifting property if and only if for any $x\in \Hom(G,k^\times)$ and any $s=\sum_{g\in G}s_gt^g\in K$ also $\rho_x(s)=\sum_{g\in G}s_gx^gt^g$ lies in $ K$.
	
	\begin{lemma}\thlabel{lem:conliftdescr2} 
		Let $x\in \Hom(G,k^\times)$ and let $r=(g_i,r_i)_{i<\alpha}\in \LRS(k,G)$. Then $\rho_x^*(r)=(g_i,x^{g_i}r_i)_{i<\alpha}$.
	\end{lemma}

	\begin{proof} 
		First note that
		$\dom(\rho_x^*(r))=\supp(\rho_x(\widetilde{r}^*))=\supp(\widetilde{r}^*)=\dom(r)$. Since $\rho_x(r^*)=\sum_{i<\alpha}r_ix^{g_i}t^{g_i}$, we obtain for any $i<\alpha$ that
		$\rho_x^*(r)(g_i)=r_ix^{g_i}$, as required.
	\end{proof}

	\begin{example} \thlabel{ex:dualitylift2}
		\begin{enumerate}[label = (\roman*)]
			\item 	Let $G=\Z$ and consider $x\in \Hom(\Z,k^\times)$ induced by $x^1=a$ for some $a\in k^\times$. Then for any $r=(g_i,r_i)_{i<\alpha}\in \LRS(k,\Z)$, we have $\rho_x^*(r)=(g_i,a^{g_i}r_i)_{i<\alpha}$.
			
			\item\label{ex:dualitylift2:2} Let $A\subseteq G$ be well-ordered
			and consider $r_A$ from \Autoref{not:rayner}.
			Then for any $x\in \Hom(G,k^\times)$, the linear recurrence sequence $\rho_x^*(r)$ has domain $A\cup \{0\}$ and is given by $\rho_x^*(r_A)(0)=x^0r_A(0)=1$ and $\rho_x^*(r_A)(g)=x^gr_A(g)=0$ for any $g\in A\setminus\{0\}$. Hence, $\rho_x^*(r_A)=r_{A}$.
		\end{enumerate}
		\qed
	\end{example}

	\begin{theorem}\thlabel{thm:2lpcriterion} 
		Let $R\subseteq \LRS(k,G)$ such that $K=\gen{R}$ is a Hahn field. Suppose that for any $x\in \Hom(G,k^\times)$ we have $\rho_x^*(R)\subseteq R$.
		Then $K$ has the canonical second lifting property. 
	\end{theorem}

	\begin{proof} 
		Let $x\in \Hom(G,k^\times)$.
		Since $\supp(\rho_x(a))=\supp(a)$ for any $a\in k\pow{G}$, we can apply \Autoref{prop:phi*lift} to obtain 
		$\rho_x(\gen{r})=\gen{\rho_x^*(r)}$. This already implies  $\rho_x(K)\subseteq K$, as required.
	\end{proof}

	We now examine how \Autoref{thm:2lpcriterion} can be applied to Hahn fields determined by $\cF$-sequences as well as to $k$-hulls.

	\begin{corollary}\thlabel{thm:canonicallifting2} 
		\begin{enumerate}[label = (\roman*)]
			\item\label{thm:canonicallifting2:1} If $\gen{S(\cF)}$ is a Hahn field, then it has the canonical second lifting property.
			\item\label{thm:canonicallifting2:2} If $\cF_{\fin}\subseteq \cF$ and $\cF$ is closed under addition, then $\gen{S(\cF)}$ is a Hahn field with the canonical second lifting property.
		\end{enumerate}
	\end{corollary}
	
	\begin{proof} 
		\Autoref{cor:hahnfield} shows that \autoref{thm:canonicallifting2:1} implies \autoref{thm:canonicallifting2:2}. In order to prove \autoref{thm:canonicallifting2:1},  let $x\in \Hom(G,k^\times)$. Moreover, let $A\in \cF$ and let $r=(g_i,r_i)_{i<\alpha}\in \LRS(k,G)$ be non-trivial with $\dom(r)=A$.
		Then  $\rho_x^*(r)=(g_i,x^{g_i}r_i)_{i<\alpha}$ by \Autoref{lem:conliftdescr2}. Thus, $\rho_x^*(r)$ is non-trivial and its domain is given by $A\in \cF$, i.e.\ $\rho_x^*(r)$ is an $\cF$-sequence, as required
		by \Autoref{thm:2lpcriterion}.
	\end{proof}

	Note that \Autoref{thm:canonicallifting2}~\autoref{thm:canonicallifting2:1} can be applied to $k\pow{\x}$ for any $n\in \N\setminus\{0\}$: since $k\pow{\x}=\gen{S(\{\N^n\})}$ and $k\pow{\x}$ is a Hahn field by \Autoref{rmk:k((x))props}, it has the canonical second lifting property.

	\begin{corollary}\thlabel{cor:canonical2lifting2} 
		Let $k\neq \mathbb{F}_2$ and suppose that $k\pow{\cF}$ is a field. Then it is a Hahn field with the canonical second lifting property if and only if $\{\{g\}\mid g\in G\}\subseteq \cF$.
	\end{corollary}
	
	\begin{proof} 
		
		Arguing as in the proof of \Autoref{cor:canonicallifting2}, we have that $\gen{R_{\cF}}=k\pow{\cF}$ and that
		$k\pow{\cF}$ is a Hahn field if and only if $\{\{g\}\mid g\in G\}\subseteq \cF$.
		
		Let $x\in \Hom(G,k^\times)$ and let $A\in \cF$. Then by
		\Autoref{ex:dualitylift2}~\autoref{ex:dualitylift2:2}, we have $\rho_x^*(r_A)=r_{A}\in R_{\cF}$. Hence, $k\pow{\cF}$ has the canonical second lifting property by \Autoref{thm:2lpcriterion}.
	\end{proof}

	\begin{remark} 
		\Autoref{cor:canonical2lifting2} was proved for Rayner fields (also with $k=\mathbb{F}_2$) in \cite[Proposition~3.6.5]{kuhlmann-serra-fields}, that is, any Rayner field has the canonical second lifting property. Thus by
		\Autoref{thm:canonicallifting2}, all Hahn fields determined by linear recurrence relations from 
		\Autoref{constr:1clpfields} also have the canonical second lifting property.\qed
	\end{remark}

	

	\begin{footnotesize}
		
	\end{footnotesize}


\begin{thebibliography}{99}			
			
			\bibitem{ax}
			\textsc{J.~Ax},
			`On Schanuel's conjectures', 
			\textsl{ Ann.\ of Math.\ (2)} 93 (1971) 252--268,
			doi:10.2307/1970774.
			
			\bibitem{conrad}
			\textsc{P.~Conrad},
			`The group of order preserving automorphisms of an ordered abelian group', 
			\textsl{Proc.\ Amer.\ Math.\ Soc.} 9 (1958) 382--389,
			doi:10.2307/2032993.
			
			\bibitem{duchamp}
			\textsc{G.~Duchamp} and \textsc{C.~Reutenauer},
			`Un critère de rationalité pro\-ve\-nant de la géométrie non commutative',
			\textsl{Invent.\ Math.} 128 (1997) 613--622,
			doi:10.1007/s002220050154.
				
			\bibitem{hahn}
			\textsc{H.~Hahn},
			`Über die nichtarchimedischen Grössensysteme', 
			\textsl{Wien.\ Ber.} 116 (1907) 601--655.
			
			\bibitem{hardy}
			\textsc{G.~H.~Hardy},
			\textsl{A course of pure mathematics},
			10th edn (Cambridge Univ.\ Press, London, 1952).
			
			\bibitem{hickel}
			\textsc{M.~Hickel} and \textsc{M.~Matusinski},
			`On the algebraicity of Puiseux series',
			\textsl{Rev.\ Mat.\ Complut.} 30 (2017) 589--620,
			doi:10.1007/s13163-017-0236-3.
			
			\bibitem{hickel2}
			\textsc{M.~Hickel} and \textsc{M.~Matusinski},
			`About algebraic Puiseux series in several variables',
			\textsl{J.\ Algebra} 527 (2019) 55--108,
			doi:10.1016/j.jalgebra.2019.02.004.
			
			\bibitem{hickel3}
			{\textsc{M.~Hickel} and \textsc{M.~Matusinski},
			`About the algebraic closure of formal power series in several variables',
			Preprint, 2023,
			arXiv:2307.04424v1.}
			
			\bibitem{kks}
			\textsc{L.~S.~Krapp}, \textsc{S.~Kuhlmann} and \textsc{M.~Serra},
			`On Rayner structures',
			\textsl{Comm.\ Algebra} 50 (2022) 940--948, 
			doi:10.1080/00927872.2021.1976789.
			
			\bibitem{kronecker}
			\textsc{L. Kronecker}, `Zur Theorie der Elimination einer Variabeln aus zwei algebraischen
			Gleichungen', \textsl{Berl. Monatsber.} (1881), 535--600.
			
			\bibitem{kuhlmann}
			\textsc{S.~Kuhlmann} and \textsc{S.~Shelah},
			`$\kappa$-bounded exponential-logarithmic power series fields',
			\textsl{Ann.\ Pure Appl.\ Logic} 136 (2005) 284--296,
			doi:10.1016/j.apal.2005.04.001.
			
			\bibitem{kuhlmannmatusinskipoint}
			\textsc{S.~Kuhlmann}, \textsc{M.~Matusinski} and \textsc{F.~Point},
			`The Valuation Difference Rank of a Quasi-Ordered Difference Field',
			\textsl{Groups, Modules, and Model Theory -- Surveys and Recent
				Developments} (eds M.\ Droste, L.\ Fuchs, B.\ Goldsmith and
			L.\ Strüngmann; Springer, Cham, 2017) 399--414,
			doi:10.1007/978-3-319-51718-6\_23.
			
			\bibitem{kuhlmann-serra-fields}
			\textsc{S.~Kuhlmann} and \textsc{M.~Serra},
			`The automorphism group of a valued field of generalised power series',
			\textsl{J. Algebra} 605 (2022) 339--376, doi:10.1016/j.jalgebra.2022.04.023.
			
			\bibitem{mourrain}
			\textsc{B.~Mourrain},
			`Fast Algorithm for Border Bases of Artinian Gorenstein Algebras',
			\textsl{ISSAC'17: Proceedings of the 2017 ACM International Symposium on Symbolic and Algebraic Computation}, Kaiserslautern, 2017 (ed.\ M.~Burr; ACM, New York, 2017) 333--340,
			doi:10.1145/3087604.3087632.
						
			\bibitem{neumann} 
			\textsc{B.~H.~Neumann},
			`On ordered division rings',
			\textsl{Trans.\ Amer.\ Math.\ Soc.} 66 (1949), 202--252, 
			doi:10.1090/S0002-9947-1949-0032593-5.
			
			\bibitem{power}
			\textsc{S.~C.~Power},
			`Finite Rank Multivariable Hankel Forms',
			\textsl{Linear Algebra Appl.} 48 (1982) 237--244,
			doi:10.1016/0024-3795(82)90110-0.

			\bibitem{reutenauer}
			\textsc{C.~Reutenauer},
			`Michel Fliess and non-commutative formal power series',
			\textsl{Int.\ J.\ Control} 81 (2008) 336--341,
			doi:10.1080/00207170701556898.

			\bibitem{ribenboim}
			\textsc{P.~Ribenboim},
			\textsl{The Theory of Classical Valuations},
			Springer Monogr.\ Math. (Springer, New York, 1999).
			
			\bibitem{salem}
			\textsc{R.~Salem},
			\textsl{Algebraic Numbers and Fourier Analysis},
			Heath Math.\ Monogr.\ (Heath, Boston, 1963).
			
			\bibitem{serra-thesis}
			\textsc{M.~Serra},
			`Automorphism groups of Hahn groups and Hahn fields', Doctoral Thesis, Universit{\"a}t Konstanz, 2021.
		\end{thebibliography}
\end{document}

pdflatex: --aux-directory=build
bibtex   : build/